\documentclass[twoside,notitlepage,11pt]{article}

\usepackage{amssymb}
\usepackage[leqno]{amsmath}
\usepackage{amsfonts}
\usepackage{amsopn}
\usepackage{amstext}
\usepackage{amsthm}
\usepackage{makeidx}
\usepackage{verbatim}
\usepackage[colorlinks, backref]{hyperref}

\newcommand{\sm}{\setminus}

\providecommand{\cal}{\mathcal}
\renewcommand{\Bbb}{\mathbb}

\newenvironment{pf}{\begin{proof}}{\end{proof}}

\newcommand{\Ef}{{\cal{F}}}

\newcommand{\Nat}{{\Bbb{N}}}

\newcommand{\al}{\alpha}

\newcommand{\eps}{\varepsilon}
\renewcommand{\phi}{\varphi}
\renewcommand{\rho}{\varrho}

\newcommand{\rest}{\restriction}

\newcommand{\loe}{\leq}
\newcommand{\goe}{\geq}

\newcommand{\subs}{\subseteq}
\newcommand{\sups}{\supseteq}
\newcommand{\nnempty}{\ne\emptyset}
\newcommand{\argum}{\:\cdot\:}
\newcommand{\ovr}{\overline}

\newcommand{\nbd}{\operatorname{nbd}}

\newcommand{\id}[1]{{\operatorname{i\!d}_{#1}}} % identity morphism

\newcommand{\dom}{\operatorname{dom}}

\newcommand{\rng}{\operatorname{rng}}
\newcommand{\supp}{\operatorname{supp}}

\newcommand{\oraz}{\qquad\text{and}\qquad}

\newcommand{\length}{\operatorname{length}}

\newtheorem{tw}{Theorem}[section]
\newtheorem{wn}[tw]{Corollary}
\newtheorem{lm}[tw]{Lemma}
\newtheorem{prop}[tw]{Proposition}
\newtheorem{claim}[tw]{Claim}
\newtheorem{subclaim}[tw]{Subclaim}
\theoremstyle{definition}
\newtheorem{df}[tw]{Definition}

\newtheorem{uwgi}[tw]{Remark}

\newcommand{\setof}[2]{\{#1\colon #2\}}
\newcommand{\setdots}[2]{\{#1,\dots,#2\}}
\newcommand{\bigsetof}[2]{\Bigl\{#1\colon #2\Bigr\}}

\newcommand{\sett}[2]{\{#1\}_{#2}}
\newcommand{\sn}[1]{\{#1\}} % singleton
\newcommand{\dn}[2]{\{#1,#2\}} % doubleton
\newcommand{\pair}[2]{\langle #1, #2 \rangle} % pair
\newcommand{\map}[3]{#1\colon #2 \to #3} % A function
\newcommand{\img}[2]{#1[#2]} % image of a set
\newcommand{\U}{\mathbb U}

\newcommand{\bbU}{{\ensuremath{\mathbb{U}}}}
\newcommand{\bbR}{{\ensuremath{\mathbb{R}}}}
\newcommand{\mc}{{\ensuremath{\operatorname{MC}}}}
\newcommand{\lip}{\operatorname{LIP}}
\newcommand{\iso}[3]{#1\colon #2 \cong #3} % An iso

\newcommand{\cmp}{\circ}
\newcommand{\num}[1]{$(#1)$}

\newcommand{\sLip}{\ensuremath{\Gamma^{\operatorname{LIP}}}}
\newcommand{\sHld}{\ensuremath{\Gamma^{\operatorname{HLD}}}}
\newcommand{\homeo}{\operatorname{H}}
\newcommand{\homeolc}[1]{\operatorname{H}_{{#1}}^{\operatorname{LC}}}
\newcommand{\lipc}{\operatorname{lip}}

\newcommand{\separator}{\begin{center}***\end{center}}
\newcommand{\abs}[1]{|#1|}
\newcommand{\calA}{{\mathcal A}}
\newcommand{\almostcontained}{\subs^*}
\newcommand{\stablep}[2]{\operatorname{S}({#1},{#2})}
\newcommand{\til}{\tilde}

\newcommand{\lhs}[1]{{\operatorname{lhs}(\mbox{#1})}}
\newcommand{\rhs}[1]{{\operatorname{rhs}(\mbox{#1})}}

\title{Extension and reconstruction theorems for the Urysohn
universal metric space}
\author{
{Wies{\l}aw Kubi\'s}\footnote{Corresponding author. E-mail address: \texttt{kubis@math.cas.cz}}\\
{\small Institute of Mathematics}\\
{\small Czech Academy of Sciences}\\
{\small Prague, CZECH REPUBLIC}\\
{\small {\it and }}\\
{\small Institute of Mathematics}\\
{\small Jan Kochanowski University in Kielce}\\
{\small POLAND}
\and
{Matatyahu Rubin}\\
{\small Department of Mathematics}\\
{\small Ben Gurion University of the Negev}\\
{\small Be'er Sheva, ISRAEL}
}

\begin{document}

\maketitle

\begin{abstract}
%We prove an extension theorem involving bilipschitz maps of open balls in the Urysohn space. We also prove an extension theorem for maps which are continuous with respect to a fixed modulus of continuity. Our results are then applied to reconstruction theorems of groups of locally bilipschitz homeomorphisms of open subsets of the Urysohn space.
We prove some extension theorems involving uniformly continuous maps of the universal Urysohn space. As an application, we prove reconstruction theorems for certain groups of autohomeomorphisms of this space and of its open subsets.
\\

{\bf Keywords and phrases:} Urysohn space, bilipschitz homeomorphism, modulus of continuity, reconstruction theorem, extension theorem.
\\

{\bf Mathematics Subject Classification (2000):}
22F50, % Noncompact transformation groups -- Groups as automorphisms of other structures
54E40, % Special maps on metric spaces
51F99, % Metric geometry -- None of the above, but in this section (???)
20E36, % Group theory -- General theorems concerning automorphisms of groups
54H11. % Topological groups 
\end{abstract}

\tableofcontents

\section{Introduction}

This work deals with the {\it Urysohn space} \cite{Ur}, which we denote by \bbU.
\index{\bbU}\index{Urysohn space}
This is the unique, up to isometry, complete separable metric space with the following properties.
\begin{enumerate}
	\item[(1)] Every separable metric space is isometrically embeddable in $\bbU$.
	\item[(2)] For every $A,B \subs \bbU$ and $\map{f}{A}{B}$:
if $A$ is finite and $f$ is an isometry between $A$ and $B$,
then there is an isometric bijection $g$ from $\bbU$ to $\bbU$
such that~$f \subs g$.
\end{enumerate}
\index{$\lip(\argum)$}
\index{$\homeo(\argum)$}
We investigate the group $\lip(\bbU)$
of bilipschitz homeomorphisms of $\bbU$ and some related groups.
Indeed, the group $\homeo(\bbU)$ of all homeomorphisms of $\bbU$
comes to mind first.
However, by the result of V. Uspenskiy \cite{Us}, $\bbU$ is homeomorphic to $\ell_2$.
So $\homeo(\bbU)$ is in fact the group of homeomomorphisms
of a Banach space, and can be better understood as such.
For $\lip(\bbU)$ and for other groups defined via the metric
of $\bbU$, the fact that $\bbU \cong \ell_2$ does not seem to help.

The main tool and also the main result in this work
is an extension theorem for finite bilipschitz functions
defined on subsets of $\bbU$ (Theorem~\ref{t2.1}).
Suppose that $A$ is a finite subset of $\bbU$ and
$\map{f}{A}{\bbU}$ is $K$-bilipschitz.
We shall show that there is $g \in \lip(\bbU)$
such that $g \sups f$, $g$ is $K$-bilipschitz
and for some open ball $B$, $g \rest (\bbU\sm B) = \id{}$, where $\id{}$ denotes the identity map.

We now turn to the description of the other groups considered in this
work. Our interest in these groups is two-fold: extension theorems
and reconstruction theorems.
By a reconstruction theorem we mean a statement of the following
form: If $\iso{\varphi}{\lip(X)}{\lip(Y)}$,
that is, if $\varphi$ is an isomorphism between the groups
$\lip(X)$ and $\lip(Y)$,
then there is a bilipschitz homeomorphism $\tau$ between $X$ and $Y$
such that for every $g \in \lip(X)$,
$\varphi(g) = \tau \cmp g \cmp \tau^{-1}$.

Recall that a {\it modulus of continuity} is a concave homeomorphism of
$[0,\infty)$, i.e. a homeomorphism $\alpha$ of $[0,\infty)$ satisfying
$\alpha(\lambda x + (1 - \lambda)y) \geq
\lambda \alpha(x) + (1 - \lambda) \alpha(y)$
for every $x,y \in [0,\infty)$ and $\lambda \in (0,1)$.
Denote by $\mc$ the set of all moduli of continuity.
\index{modulus of continuity}
\index{\mc}
\index{$\al$-continuous function}
Given $\al\in\mc$, we say that a function $f$ from a metric space $(X,d^X)$
to a metric space $(Y,d^Y)$ is {\it $\al$-continuous}
if $d^Y(f(u),f(v)) \leq \alpha(d^X(u,v))$
for every $u,v \in X$. 

We do not know to generalize the extension theorem for finite
bilipschitz functions (Theorem~\ref{t2.1}),
to a general $\alpha \in \mc$.
Whereas for bilipschitz functions we prove that every
finite $K$-bilipschitz function can be extended
to a $K$-bilipschitz homeomorphism of $\bbU$
which is the identity outside a ball,
for a general $\alpha \in \mc$, we only know to prove
that every finite $\alpha$-bicontinuous function is extendible
to an $\alpha$-bicontinuous homeomorphism of $\bbU$
(Corollary~\ref{c3.3}).

The fact that in bilipschitz case the extending homeomorphism
can be constructed in such a way
that it is the identity outside a ball, means that we also get
an ``Extension theorem for finite bilipschitz functions''
for open subsets of $\bbU$ and not just for $\bbU$.

However, both Theorem~\ref{t2.1} and Corollary~\ref{c3.3}
can be strengthened by proving the extension theorem
not just for functions with a finite domain,
but also for functions whose domain is a totally bounded set.
For isometries, this fact is due to Huhunai\v svili \cite{Hu}.
The argument we use is similar to Huhunai\v svili's.
These results appear in Theorems~\ref{t2.91} and \ref{t3.4}.

In order to state the next result, let us give necessary definitions.

\index{$\preceq$}
Fix $\al,\beta\in\mc$. We shall write $\al\preceq \beta$ if there is $a>0$ such that $\al\rest[0,a]\leq \beta\rest[0,a]$, where $f\leq g$ means that $f(x)\leq g(x)$ for every $x$ in the common domain of $f$, $g$. Note that when $K>0$ then the function $x\mapsto Kx$ belongs to \mc. Also, if $\al,\beta\in \mc$ then $\al+\beta\in \mc$ and $\al\cmp \beta\in \mc$.

\index{measure of continuity semigroup}\index{\mc-semigroup}
A subset $\Gamma$ of \mc\ will be called a {\it measure of continuity semigroup} (briefly: {\it \mc-semigroup}) if the following holds.
\begin{enumerate}
	\item[(1)] $\Gamma$ contains the function $y = 2x$.
	\item[(2)] $\Gamma$ is closed under compositions, i.e. $\al\cmp\beta\in \Gamma$ whenever $\al,\beta\in\Gamma$.
	\item[(3)] For every $\al\in\Gamma$, $\setof{\beta\in\mc}{\beta\preceq \al}\subs\Gamma$.
\end{enumerate}
\index{\mc-semigroup!-- countably generated}
Further, we say that $\Gamma$ is {\it countably generated} if there exists a countable set $\Gamma_0\subs\Gamma$ such that for every $\al\in\Gamma$ there is $\beta\in\Gamma_0$ with $\al\preceq \beta$. In other words, $\Gamma = \setof{\al\in\mc}{(\exists\;\beta\in\Gamma_0)\;\al\preceq\beta}$.

\index{\mc-semigroup!-- Lipschitz}\index{\mc-semigroup!-- H\"older}\index{\sLip}\index{\sHld}
Important examples of countably generated \mc-semigroups are {\it Lipschitz \mc-semigroup} \sLip\ and {\it H\"older semigroup} \sHld. The first one is generated by functions of type $x\mapsto nx$ ($n\in\Nat$) and the latter one is generated by functions of the form $x\mapsto x^{1/n}$ ($n\in\Nat$). Note that every \mc-semigroup contains \sLip.

\index{$\Gamma$-continuous map}\index{locally $\Gamma$-bicontinuous map}
Let $\Gamma$ be an MC-semigroup and $f$ be a
function from a metric space $X$ to a metric space $Y$.
Then $f$ is {\it locally $\Gamma$-continuous}
if for every $x \in X$ there is a neighborhood $U$ of $x$
and $\alpha \in \Gamma$ such that $f \rest U$
is $\alpha$-continuous.
The function $f$ is {\it locally $\Gamma$-bicontinuous},
if $f$ is a homeomorphism between
$X$ and $\rng(f)$, and both $f$ and $f^{-1}$ are locally
$\Gamma$-continuous.

From the extension theorem for finite bilipschitz functions
we shall deduce the following result (Corollary~\ref{c6.4}(a)).

\begin{tw}\label{t1.4}
Let $X$ and $Y$ be open subsets of $\bbU$,
$\Gamma$ and $\Delta$ be countably generated MC-semigroups and
$\iso{\varphi}
{\homeolc{\Gamma}(X)}{\homeolc{\Delta}(Y)}$.
Then there is $\iso{\tau}{X}{Y}$
such that $\tau$ is locally $\Gamma$-bicontinuous,
$\tau$ is locally $\Delta$-bicontinuous
and $\varphi(g) = \tau \cmp g \cmp \tau^{-1}$
for every $g \in \homeolc{\Gamma}(X)$.
\end{tw}

Note that the above theorem does not claim that $\Gamma = \Delta$.
However, the theorem does imply the weaker statement that
$\homeolc{\Gamma}(X) = \homeolc{\Delta}(X)$.
We know to conclude that, indeed, $\Gamma = \Delta$,
only when $X = Y = \bbU$.
It is Theorem~\ref{t3.4} which is used in order to deduce this.

In fact, we do not know to prove
the following general statement.
\begin{quote}
If $X$ is an open subset of $\bbU$
and $\Gamma$ and $\Delta$ are different MC-semigroups,
then $\homeolc{\Gamma}(X) \neq \homeolc{\Gamma}(Y)$.
\vspace{-05.7pt}
\end{quote}
Even the following special case is unknown.
For $\alpha \in \mc$ let $\Gamma_{\alpha}$
denote the $\mc$-semigroup generated by $\alpha$.
Suppose that $\Delta$ is a countably generated $\mc$-semigroup
and $\alpha \not\in \Delta$.
Is it true that
$\homeolc{\Gamma_{\alpha}}(X) \not\subs
\homeolc{\Delta}(X)$?
Note that if $\bbU$ is replaced by $\bbR$,
that is, we consider an open subset of $\bbR$,
then different MC-semigroups do define different subgroups
of $\homeo(X)$.
The same is true for any open subset of a normed linear space.

%End of Introduction

\section{Extensions of finite bilipschitz functions}

The main result of this section is Theorem \ref{t2.1}.
It is the basis for our reconstruction theorem for open subsets
of $\bbU$.

\index{$\id{}$}
Let $B(x,r) \subs \bbU$ be an open ball in $\bbU$.
We wish to prove the following claim.
Let $A \subs B(x,r)$ be a finite set and $K > 1$.
Assume that $\map fA{B(x,r)}$ is such that $f \cup \id{\bbU \sm B(x,r)}$ is
$K$-bilipschitz, where $\id S$ denotes the identity map on the set $S$.
Then there is $g \in \homeo(\bbU)$ such that
$g \sups f \cup \id{\bbU \sm B(x,r)}$
and $g$ is $K$-bilipschitz. It turns out that we need some extra assumptions in order to prove such a claim. These assumptions are somewhat technical, although good enough for applications.

\index{$N$-good function}\index{$N$-bigood function}
Fix $r > 0$ and let $f$ be a function whose domain and range are
subsets of $B(x,r)$. Also, let $N > 1$.
We shall say that $f$ is {\it $N$-good}, if
$$d(y,f(y)) \leq
\frac{1}{N} \cdot (r - d(y,x))
\mbox{ \,for every $y \in \dom(f)$.}
$$
The function $f$ is {\it $N$-bigood},
if both $f$ and $f^{-1}$ are $N$-good.

\begin{tw}[Bilipschitz Extension Theorem for the Urysohn space]\label{t2.1}
Let $B(x,r)$ be an open ball in $\bbU$,
$N \geq 4$,
$\frac{N - \sqrt{N^2 - 4N}}{2} \leq K \leq
\frac{N + \sqrt{N^2 - 4N}}{2}$,
$A \subs B(x,r)$ be finite and $x \in A$.
Suppose that $\map{f}{A}{B(x,r)}$, $f$ is $K$-bilipschitz,
$f$ is $N$-bigood and $f(x) = x$.

Then there exists $\map{g}{B(x,r)}{B(x,r)}$ such that
\begin{itemize}
\item[\num{1}] $g$ is a bijection and $g \sups f$,
\item[\num{2}] $g \cup \id{\bbU \sm B(x,r)}$ is $K$-bilipschitz,
\item[\num{3}] $g$ is $N$-bigood.
\end{itemize}
\end{tw}

Let us first see that the assumptions in the above statement are meaningful. Note that
\begin{equation}
\lim_{N \rightarrow \infty} \frac{N - \sqrt{N^2 - 4N}}{2} = 1\oraz \lim_{N \rightarrow \infty} \frac{N + \sqrt{N^2 - 4N}}{2} = \infty.
\notag
\end{equation}
It follows that for every $K>1$ there is $N\geq4$ such that $K,N$ fulfill the assumptions in Theorem \ref{t2.1}. Thus, Theorem~\ref{t2.1} can be applied to any $K > 1$
provided that the distance between any point $u \in A$
and its image $f(u)$ is sufficiently small in comparison with the
distance of $u$ and $f(u)$ from the boundary of $B(x,r)$.

The following simple facts related to the assumptions of the above theorem will be used later.

\begin{prop}\label{dwadwa}
(a) Assume $K>1$ and $N\geq \frac{K^2}{K-1}$. Then $N\geq4$ and 
\begin{equation}
\frac{N - \sqrt{N^2 - 4N}}{2}\leq K\leq \frac{N + \sqrt{N^2 - 4N}}{2}.
\tag{*}\label{weiothjqwer}
\end{equation}

(b) If $N\geq4$ and (\ref{weiothjqwer}) holds then $K>1$ and $N\geq\frac{K^2}{K-1}$.

(c) $\frac1N + \frac1K \leq1$ iff $1 + \frac1{N-1} \leq K$.

(d) The inequalities in (c) hold whenever $K,N$ satisfy the requirements of Theorem \ref{t2.1}.
\end{prop}

\index{compliant function}
A function $f$ whose domain and range are subsets of $B(x,r)$ will be called briefly
{\it $(K,N)$-compliant} if it is both $K$-bilipschitz and $N$-bigood.
Note that $f$ is $(K,N)$-compliant iff
$f^{-1}$ is $(K,N)$-compliant.

\begin{prop}\label{p2.2}
Let $X$ be a metric space, $x\in X$, $r>0$.
Suppose that $h$ is a function whose domain and range are subsets of $B(x,r)$.
Let $K,N  > 0$ and assume that $1 + \frac{1}{N} \leq K$.
If $h$ is $(K,N)$-compliant then $h\cup\id{X\sm B(x,r)}$ is $K$-bilipschitz.
\end{prop}

\begin{pf}
Fix $u \in \dom(h)$ and fix $w \in X \sm B(x,r)$.
Then $r - d(x,u) \leq d(w,x) - d(x,u) \leq d(w,u)$.
Hence $r - d(x,u) \leq d(u,w)$. Thus
\begin{align*}
d(h(u),h(w)) &= d(h(u),w) \leq d(h(u),u) + d(u,w) \leq
\frac{1}{N} \cdot (r - d(x,u)) + d(u,w) \\ &\leq
\frac{1}{N} \cdot d(u,w) + d(u,w) = (1 + \frac{1}{N}) \cdot d(u,w) \leq
K d(u,w).
\end{align*} 
\end{pf}

\index{amalgamation}
Before we prove Theorem \ref{t2.1}, we recall an important property of (isometric embeddings of) metric spaces, called {\it amalgamation}. Namely, given finite metric spaces $\pair{X_0}{d_0}$ and $\pair{X_1}{d_1}$ such that both metrics coincide on the intersection $Z=X_0\cap X_1$, there exist a metric space $(Y,d)$ and isometric embeddings $\map{f_i}{\pair{X_i}{d_i}}{\pair Yd}$, $i=0,1$, such that $f_0\rest Z=f_1\rest Z$. In other words, there is a metric on $X:=X_0\cup X_1$ which extends both $d_0$ and $d_1$. The space $\pair Xd$ will be called the {\it amalgamation} of $\pair{X_0}{d_0}$ and $\pair{X_1}{d_1}$. The amalgamation property is essentially used for constructing the Urysohn space.
We shall need the following more specific statement, which actually gives inductive argument for amalgamating two finite metric spaces.

\begin{prop}\label{p2.3ab}
Let $\pair{X_i}{d_i}$ for $i=0,1$ be finite metric spaces such that $Z=X_0\cap X_1$ is nonempty and $d_0\rest Z^2 = d_1\rest Z^2$. Assume further that $X_i=Z\cup\sn{p_i}$ for $i<2$. Then
$$\ell:=\max_{z\in Z}|d_0(p_0,z)-d_1(z,p_1)|\leq \min_{z\in Z}(d_0(p_0,z)+d_1(z,p_1))=:r$$
and every extension $d$ of $d_0\cup d_1$, satisfying 
$\ell\leq d(p_0,p_1)=d(p_1,p_0)\leq r$ and $0<d(p_0,p_1)$, is a metric on $X_0\cup X_1$.
\end{prop}

\index{amalgamation!-- minimal}
It may happen that $\ell=0$ in the above statement and then setting $d(p_0,p_1)=d(p_0,p_1)=\ell$ we also obtain amalgamation of $X_1$, $X_2$ in which points $p_0$, $p_1$ are identified. Amalgamation satisfying $d(p_0,p_1)=\ell$ will be called {\it minimal}.

The properties of the Urysohn space imply the following: given a finite metric space $X\subs \U$ and its finite metric extension $Y\sups X$, there exists an isometric embedding $\map hY\U$ such that $h\rest X=\id X$.

The crucial argument in constructing the homeomorphism $g$ promised
in Theorem~\ref{t2.1} is showing how to add a single point to
the domain or range of a $(K,N)$-compliant function so that the
resulting new function remains compliant.
Note that the problem of adding a point to the domain or to the range
of $f$ is the same, since $f$ is compliant iff $f^{-1}$ is
compliant.
%This crucial argument is the contents of the next lemma.

\begin{lm}\label{l2.4}
Assume that $N \geq 4$, $r>0$ and
$\frac{N - \sqrt{N^2 - 4N}}{2} \leq K \leq
\frac{N + \sqrt{N^2 - 4N}}{2}$.
Suppose that $f$ is such that $\dom(f),\rng(f)$
are finite subsets of $B(x_1,r)$, $f(x_1) = x_1$ and $f$ is
$(K,N)$-compliant. Let $x \in B(x_1,r)$.
Then there is $y \in B(x_1,r)$ such that
$f \cup \sn{\pair{x}{y}}$ is $(K,N)$-compliant.
\end{lm}

\begin{pf}
Let $\dom(f)=\setdots{x_1}{x_n}$, where $x_1,\ldots,x_n$ are pairwise ditstinct. Let $y_m=f(x_m)$.
Denote $d_{i,j} = d(x_i,x_j)$, $e_{i,j} = d(y_i,y_j)$ and
$s_{i,j} = d(x_i,y_j)$. The assumptions concerning $N$-bigoodness say that 
for every $m = 1,\ldots,n$, 
\begin{equation}
s_{m,m} \leq \frac{1}{N} \cdot (r - d_{m,1})\oraz s_{m,m}\leq \frac{1}{N} \cdot (r - e_{m,1}).
	\tag{B}\label{bidobroc}
\end{equation}
Assume $x \in B(x_1,r)\sm \setdots{x_1}{x_n}$ and define $d_m=d(x,x_m)$, $s_m=d(x,y_m)$. We take an imaginary point $y$ which we shall later identify, using amalgamation, with a suitable element of $B(x_1,r)$. Following is the crucial step.

\begin{claim}\label{klejm1}
There exist $e_1,\dots,e_n>0$ such that, defining $d(y,y_i)=e_i$ for $i=1,\dots,n$, the set $\{y,y_1,\dots,y_n\}$ becomes a metric space, the function $f\cup\sn{\pair xy}$ is $K$-bilipschitz and for every $m=1,\dots,n$ the following inequality holds.
\begin{equation}
\abs{e_m-s_m} \le \min\left( \frac1N (r-d_1), \frac1N (r-e_1) \right).
\tag{G}\label{eqgood}
\end{equation}
\end{claim}

Condition (\ref{eqgood}) is necessary for $N$-bigoodness of the extension: given $m\le n$, the distance between $x$ and $y$ must be at least $\abs{e_m-s_m}$ and, on the other hand, it must not exceed the right-hand side of (\ref{eqgood}).

Suppose we have proved the above claim. 

%Let $d$ be the unique two-place function on $Y'=Y\cup \sn y$, where $Y=\setdots{y_1}{y_n}$, which extends the metric of $Y$ and satisfies $d(y_m,y)=.$

We amalgamate the two metric spaces $X'=Y\cup \sn x$ and $Y'=Y\cup \sn y$, where $Y=\setdots{y_1}{y_n}$ and the metric $d$ on $Y'$ is given by Claim \ref{klejm1}, that is, $d$ extends the metric of $Y$ and $d(y_m,y)=e_m$ for $m\le n$. Let
$$s = \min \left( \min_{1\le i\le n} (e_i + s_i), \frac1N (r-d_1), \frac1N (r-e_1) \right).$$
By (\ref{eqgood}), $\max_{1\le i\le n} \abs{e_i - s_i} \le s$. Clearly, $s \le \min_{1\le i\le n}(e_i+s_i)$, therefore by Proposition \ref{p2.3ab}, there exists a metric on $X'\cup Y'$  (which we still denote by $d$) that extends the metrics of $X'$ and $Y'$ and satisfies $d(x,y)=s$. Now, amalgamate $\U$ with the finite space $X'\cup Y'$. By this way, we may assume that $y\in \U$. In fact $y\in B(x_1,r)$, because $d(x_1,y)\le d(x_1,x)+d(x,y)=d_1+s\le d_1+\frac1N(r-d_1)<r$.
Finally, $f\cup\sn{\pair xy}$ is $(K,N)$-compliant.

It remains to prove Claim \ref{klejm1}.

\separator

The distances $e_1,\ldots,e_n$ will be defined by induction on
$m = 1,\ldots,n$.
%Before starting the induction we adopt ``temporary distances''
%$e_1' = K d_1,\ldots,e_n' = K d_n$.
%In the $m$'th inductive step $e'_m$ will be replaced by the permanent distance $e_m$.
Fix $m\leq n$ and suppose that $e_1,\ldots,e_{m - 1}$ have been defined.
We assume by induction that for every $\ell < i < m$ and $j \geq m$ the following inequalities hold.
\begin{align}
e_{\ell,j} - K d_j &\leq e_{\ell} \leq e_{\ell,j} + K d_j 
	\tag{IH1}\label{ih1}\\
\abs{e_{i,\ell} - e_{\ell}} &\leq e_i \leq e_{i,\ell} + e_{\ell}
	\tag{IH2}\label{ih2}\\
\frac{1}{K} \cdot d_{\ell} &\leq e_{\ell} \leq K d_{\ell}
	\tag{IH3}\label{ih3}\\
\abs{e_{\ell} - s_{\ell}} &\leq \frac{1}{N} \cdot (r - d_1)
	\tag{IH4}\label{ih4}\\
\abs{e_{\ell} - s_{\ell}} &\leq \frac{1}{N} \cdot (r - e_1)
	\tag{IH5}\label{ih5}
\end{align}
Condition (\ref{ih1}) consists of two of the three triangle inequalities in the
triangle whose vertices are $y_{\ell},y_j$ and the future point $y$. One may think of
$Kd_j$ as a ``temporary" distance between $y$ and $y_j$.
The inequality $Kd_j \leq e_{\ell,j} + e_{\ell}$ is not assumed, since $Kd_j$ will be replaced by a smaller final distance.
(\ref{ih2}) consists of the three inequalities in the
triangle whose vertices are $y_{\ell},y_i$ and the future point $y$.
(\ref{ih3}) is the bilipschitz condition for the pairs $x_{\ell},x$ and $y_{\ell},y$.
Finally, (\ref{ih4}) and (\ref{ih5}) are just condition (\ref{eqgood}) for $\ell<m$.
It is clear that $e_1,\dots, e_n$ constructed by this inductive procedure fulfill the requirements of Claim \ref{klejm1}.

Consider the following system of inequalities in the unknown $e_m$.
\begin{align}
e_{m,j} - K d_j \leq e_m & \leq e_{m,j} + K d_j \qquad (j>m)
	\tag{IE1$_m$}\label{ie1}\\
\abs{e_{m,\ell} - e_{\ell}} \leq e_m & \leq e_{m,\ell} + e_{\ell} \qquad (\ell<m)
	\tag{IE2$_m$}\label{ie2}\\
\frac{1}{K} \cdot d_m \leq e_m & \leq K d_m
	\tag{IE3$_m$}\label{ie3}\\
s_m - \frac{1}{N} \cdot (r - d_1) \leq e_m & \leq s_m + \frac{1}{N} \cdot (r - d_1)
	\tag{IE4$_m$}\label{ie4}\\
s_m - \frac{1}{N} \cdot (r - e_1) \leq e_m & \leq s_m + \frac{1}{N} \cdot (r - e_1)
	\tag{IE5$_m$}\label{ie5}\\
e_1 & \leq N \cdot (s_i - \frac{d_i}{K}) + r \qquad (i>1)
	\tag{IE6$_1$}\label{ie6}\\
e_1 & \leq N \cdot (K d_i - s_i) + r \qquad (i>1)
	\tag{IE7$_1$}\label{ie7}
\end{align}
Observe that a solution $e_m$ to the above system satisfies inequalities (\ref{ih1}) -- (\ref{ih5}) with $\ell=m$, i.e. the inductive step can be accomplished. It remains to show that the above system is solvable.

Given an inequality with label (e), we shall denote by $\lhs e$ and $\rhs e$ the expression on the left-hand side and on the right-hand side, respectively.

Inequalities (\ref{ie6}), (\ref{ie7}) are required for solving some of the
inequalities with $m > 1$. Namely, inequality (\ref{ie6}) is needed in the proof
of $\lhs{\ref{ie3}} \leq \rhs{\ref{ie5}}$.
Inequality (\ref{ie7}) is needed for the proof of
$\lhs{\ref{ie1}} \leq\rhs{\ref{ie5}}$,
$\lhs{\ref{ie5}} \leq \rhs{\ref{ie1}}$ and
$\lhs{\ref{ie5}} \leq \rhs{\ref{ie3}}$.

It has to be shown that each expression appearing in the left is
$\leq$ every expression appearing on the right.

It is worthwhile to note that the original assumptions in the lemma together with hypotheses (\ref{ih4}) and (\ref{ih5}) will be used only in proving inequalities involving (\ref{ie4}) and (\ref{ie5}).
The verification that the inequalities arising from
(\ref{ie1}) -- (\ref{ie4}) hold does not depend on $m$.
For inequalities involving (\ref{ie5}) we need to distinguish between the cases $m>1$ and $m=1$.

\begin{subclaim}
$\lhs{\ref{ie1}}\leq \rhs{\ref{ie2}}$, i.e. $e_{m,j} - K d_j \leq e_{m,\ell} + e_{\ell}$.
\end{subclaim}

\begin{pf}
By the triangle inequality, $e_{m,j} \leq e_{m,\ell} + e_{\ell,j}$. Thus
$e_{m,j} - K d_j \leq e_{m,\ell} + e_{\ell,j} - K d_j$. Finally, the first inequality in (\ref{ih1}) says that $e_{\ell,j} - K d_j \leq e_{\ell}$.
\end{pf}

\begin{subclaim}
$\lhs{\ref{ie2}}\leq \rhs{\ref{ie1}}$, i.e.
$\abs{e_{m,\ell} - e_{\ell}} \leq e_{m,j} + K d_j$.
\end{subclaim}

\begin{pf}
Using the fact that $e_{m,\ell} \leq e_{m,j} + e_{\ell,j}$, we get
$e_{m,\ell} - e_{\ell} \leq e_{m,j} + e_{\ell,j} - e_{\ell}$.
The first inequality in (\ref{ih1}) gives $e_{\ell,j} \leq e_{\ell} + K d_j$, therefore 
$e_{m,\ell} - e_{\ell} \leq e_{m,j} + K d_j$.

On the other hand, the second inequality in (\ref{ih1}) says that $e_{\ell} \leq e_{\ell,j} + K d_j$, therefore
$e_{\ell} - e_{m,\ell} \leq e_{\ell,j} + K d_j - e_{m,\ell} \leq e_{m,j} + K d_j$.
\end{pf}

\begin{subclaim}
$\lhs{\ref{ie1}} \leq \rhs{\ref{ie3}}$, i.e. $e_{m,j} - K d_j \leq K d_m$.
\end{subclaim}

\begin{pf}
This follows from $d_{m,j} \leq d_j + d_m$ and from the fact that $f$ is $K$-Lipschitz.
\end{pf}

\begin{subclaim}
$\lhs{\ref{ie3}} \leq \rhs{\ref{ie1}}$, i.e. $\frac{1}{K} \cdot d_m \leq e_{m,j} + K d_j$.
\end{subclaim}

\begin{pf}
Since $f^{-1}$ is $K$-Lipschitz, $\frac{1}{K} \cdot d_{m,j} \leq e_{m,j}$. Thus
$\frac 1K \cdot d_m \leq \frac 1K (d_{m,j} + d_j) \leq e_{m,j} + \frac 1K \cdot d_j \leq e_{m,j} + K d_j$. The last inequality holds because $K>1$.
\end{pf}

\begin{subclaim}
$\lhs{\ref{ie2}} \leq \rhs{\ref{ie3}}$, i.e. $\abs{e_{m,\ell} - e_{\ell}} \leq K d_m$.
\end{subclaim}

\begin{pf}
This is a particular instance of (\ref{ih1}), where $j:=m$.
\end{pf}

\begin{subclaim}
$\lhs{\ref{ie3}} \leq \rhs{\ref{ie2}}$, i.e. $\frac{1}{K} \cdot d_m \leq e_{m,\ell} + e_{\ell}$.
\end{subclaim}

\begin{pf}
Notice that $\frac 1K\cdot d_{m,l}\leq e_{m,l}$, because $f$ is $K$-Lipschitz and $\frac 1K\cdot d_\ell \leq e_\ell$, by the first inequality in (\ref{ih3}). Thus
$\frac 1K\cdot d_m \leq \frac 1K\cdot d_{m,\ell} + \frac 1K\cdot d_\ell \leq e_{m,\ell} + e_\ell$.
\end{pf}

Note that until now we have neither used assumptions (\ref{bidobroc}) nor
did we use hypotheses (\ref{ih4}) and (\ref{ih5}).

\begin{subclaim}\label{qfrqwrqrqrewrpir}
$\lhs{\ref{ie1}} \leq \rhs{\ref{ie4}}$, i.e. $e_{m,j} - K d_j \leq s_m + \frac{1}{N} \cdot (r - d_1)$.
\end{subclaim}

\begin{pf}
Triangle inequalities give $e_{m,j} \leq s_m + s_j$ and $s_j \leq d_j + s_{j,j}$, so
$e_{m,j} - K d_j \leq s_m + d_j + s_{j,j} - K d_j$.
Thus, it suffices to show that $d_j + s_{j,j} \leq K d_j + \frac 1N \cdot (r - d_1).$
By the first inequality in (\ref{bidobroc}), $s_{j,j} \leq \frac 1N \cdot (r - d_{j,1})$.
So it suffices to show that
\begin{equation}
d_j + \frac 1N \cdot (r - d_{j,1}) \leq K d_j + \frac 1N \cdot (r - d_1),
	\notag
\end{equation}
which reduces to $\frac 1N \cdot (d_1 - d_{j,1}) \leq (K-1) d_j$. The last inequality holds, because $d_1 - d_{j,1} \leq d_j$, by triangle inequality, and $\frac 1N \leq K-1$, by the assumptions on $N,K$ (see Proposition \ref{dwadwa}).
\end{pf}

\begin{subclaim}
$\lhs{\ref{ie4}} \leq \rhs{\ref{ie1}}$, i.e. $s_m - \frac{1}{N} \cdot (r - d_1) \leq e_{m,j} + K d_j$.
\end{subclaim}

\begin{pf} Note that
$s_m \leq e_{m,j} + s_j$, so it suffices to show that
$s_j \leq K d_j + \frac{1}{N} \cdot (r - d_1)$.
This is the same as in the proof of Subclaim \ref{qfrqwrqrqrewrpir}.
\end{pf}

\begin{subclaim}
$\lhs{\ref{ie2}} \leq \rhs{\ref{ie4}}$, i.e. $\abs{e_{m,\ell} - e_{\ell}} \leq s_m + \frac{1}{N} \cdot (r - d_1)$.
\end{subclaim}

\begin{pf}
We have $e_{m,\ell} \leq s_m + s_{\ell}$, so 
$e_{m,\ell} - e_{\ell} \leq s_m + \abs{e_\ell - s_\ell}$.
Similarly, $e_{m,\ell} + s_m \geq s_{\ell}$ implies that 
$e_{\ell} - e_{m,\ell} \leq s_m - s_\ell + e_\ell \leq s_m + \abs{e_\ell - s_\ell}$.
Finally, $\abs{e_{\ell} - s_{\ell}} \leq \frac{1}{N} \cdot (r - d_1)$, by (\ref{ih4}).
\end{pf}

\begin{subclaim}
$\lhs{\ref{ie4}} \leq \rhs{\ref{ie2}}$, i.e. $s_m - \frac{1}{N} \cdot (r - d_1) \leq e_{m,\ell} + e_{\ell}$.
\end{subclaim}

\begin{pf}
We have $s_m \leq e_{m,\ell} + s_{\ell}$ so the above inequality follows from $s_\ell - \frac1N \cdot (r - d_1) \leq e_\ell$, which is part of (\ref{ih4}).
\end{pf}

\begin{subclaim}
$\lhs{\ref{ie3}} \leq \rhs{\ref{ie4}}$, i.e. $\frac{1}{K} \cdot d_m \leq s_m + \frac{1}{N} \cdot (r - d_1)$.
\end{subclaim}

\begin{pf}
Recall that $\frac 1K \leq 1 - \frac 1N$ (see Proposition \ref{dwadwa}). Also, $d_m \leq s_m + s_{m,m}$ and $s_{m,m}\leq \frac{1}{N} \cdot (r - d_{m,1})$, by (\ref{bidobroc}). Thus
$$\frac 1K \cdot d_m \leq d_m - \frac 1N\cdot d_m \leq s_m + s_{m,m} - \frac 1N \cdot d_m \leq s_m + \frac 1N \cdot (r - d_{m,1} - d_m).$$
Finally, $r - d_{m,1} - d_m \leq r - d_1$ holds by the triangle inequality.
\end{pf}

\begin{subclaim}
$\lhs{\ref{ie4}} \leq \rhs{\ref{ie3}}$, i.e. $s_m - \frac{1}{N} \cdot (r - d_1) \leq K d_m$.

\end{subclaim}

\begin{pf}
Using (\ref{bidobroc}) and triangle inequalities, we have
\begin{align*}
s_m &\leq d_m + s_{m,m} \leq d_m + \frac{1}{N} \cdot (r - d_{m,1}) = d_m + \frac 1N \cdot (r - d_1) + \frac 1N \cdot (d_1 - d_{m,1})\\
&\leq d_m + \frac 1N \cdot (r - d_1) + \frac 1N\cdot d_m \leq (1+\frac 1N)\cdot d_m + \frac 1N \cdot (r - d_1).
\end{align*}
Finally $(1+\frac 1N)\cdot d_m \leq K d_m$, by Proposition \ref{dwadwa}.
\end{pf}

We now deal with inequalities involving (\ref{ie5}).
Here we have to distinguish between the cases $m = 1$ and $m > 1$.
We start with the case $m > 1$.

\begin{subclaim}\label{wehqroiqhr}
$m>1\implies \lhs{\ref{ie1}} \leq \rhs{\ref{ie5}}$, i.e. $e_{m,j} - K d_j \leq s_m + \frac{1}{N} \cdot (r - e_1)$.
\end{subclaim}

\begin{pf}
Since $e_{m,j} - s_m \leq s_j$, the above inequality follows from
$s_j \leq \frac{1}{N} \cdot (r - e_1) + K d_j$.
That is, $e_1 \leq N \cdot (K d_j - s_j) + r$.
This is just (\ref{ie7}).
\end{pf}

\begin{subclaim}
$m>1\implies \lhs{\ref{ie5}} \leq \rhs{\ref{ie1}}$, i.e. $s_m - \frac{1}{N} \cdot (r - e_1) \leq e_{m,j} + K d_j$.
\end{subclaim}

\begin{pf}
Like in the previous subclaim, using inequality $s_m - e_{m,j}\leq s_j$ instead.
\end{pf}

\begin{subclaim}
$m>1 \implies \lhs{\ref{ie2}} \leq \rhs{\ref{ie5}}$, i.e. $\abs{e_{m,\ell} - e_{\ell}} \leq s_m + \frac{1}{N} \cdot (r - e_1)$.
\end{subclaim}

\begin{pf} Using inequalities $e_{m,\ell} - s_m \leq s_\ell \leq e_{m,\ell} + s_m$, we obtain that $e_{m,\ell} - e_\ell - s_m \leq s_\ell - e_\ell$ and $e_\ell - e_{m,\ell} - s_m \leq e_\ell - s_\ell$. Thus $\abs{e_{m,\ell} - e_{\ell}} - s_m \leq \abs{s_\ell - e_\ell}$.
Finally, $\abs{s_\ell - e_\ell} \leq \frac{1}{N} \cdot (r - e_1)$ is hypothesis (\ref{ih5}).
\end{pf}

\begin{subclaim}
$m>1 \implies \lhs{\ref{ie5}} \leq \rhs{\ref{ie2}}$, i.e. $s_m - \frac{1}{N} \cdot (r - e_1) \leq e_{m,\ell} + e_{\ell}$.
\end{subclaim}

\begin{pf} Note that
$s_m - e_{m,\ell} \leq s_{\ell}$, therefore the above inequality follows from (\ref{ih5}).
\end{pf}

\begin{subclaim}
$m>1 \implies \lhs{\ref{ie3}} \leq \rhs{\ref{ie5}}$, i.e. $\frac{1}{K} \cdot d_m \leq s_m + \frac{1}{N} \cdot (r - e_1)$.
\end{subclaim}

\begin{pf}
This is equivalent to
$e_1 \leq N \cdot (s_m - \frac{d_m}{K}) + r$, which is (\ref{ie6}).
\end{pf}

\begin{subclaim}
$m>1 \implies \lhs{\ref{ie5}} \leq \rhs{\ref{ie3}}$, i.e. $s_m - \frac{1}{N} \cdot (r - e_1) \leq K d_m$.
\end{subclaim}

\begin{pf}
That is, $e_1 \leq N \cdot (K d_m - s_m) + r$. This is (\ref{ie7}).
\end{pf}

\begin{subclaim}
$m>1 \implies \lhs{\ref{ie4}} \leq \rhs{\ref{ie5}}$, i.e. $s_m - \frac{1}{N} \cdot (r - d_1) \leq s_m + \frac{1}{N} \cdot (r - e_1)$.
\end{subclaim}

\begin{pf}
This follows from $d_1+e_1 \leq 2r$, since these are distances between points in the ball $B(x_1,r)$.
\end{pf}

\begin{subclaim}\label{q2rugiyi}
$m>1 \implies \lhs{\ref{ie5}} \leq \rhs{\ref{ie4}}$, i.e. $s_m - \frac{1}{N} \cdot (r - e_1) \leq s_m + \frac{1}{N} \cdot (r - d_1)$.
\end{subclaim}

\begin{pf}
The same as in the previous subclaim.
\end{pf}

We are left with the case $m=1$. Note that (\ref{ie2}) is vacuous in case $m=1$.  
Inequalities (\ref{ie5}) with $m=1$ need to be rewritten in such a way that $e_1$ occurs only in the middle. Further, notice that $s_1 = d_1$, because $x_1 = y_1$. Thus (\ref{ie4}) can be simplified. 
Summarizing, we now have to deal with the following system of inequalities, together with (\ref{ie6}), (\ref{ie7}).
\begin{align}
e_{1,j} - K d_j \leq e_1 & \leq e_{1,j} + K d_j \qquad (j>1)
	\tag{IE1$_1$}\label{ie11}\\
\frac{1}{K} \cdot d_1 \leq e_1 & \leq K d_1
	\tag{IE3$_1$}\label{ie31}\\
\frac{N+1}{N} \cdot d_1 - \frac{1}{N} \cdot r \leq e_1 & \leq \frac{N-1}N \cdot d_1 + \frac{1}{N} \cdot r
	\tag{IE4$_1$}\label{ie41}\\
\frac{N}{N - 1} \cdot d_1 - \frac{1}{N - 1} \cdot r \leq e_1 & \leq
\frac{N}{N + 1} \cdot d_1 + \frac{1}{N + 1} \cdot r
	\tag{IE5$_1$}\label{ie51}\\
e_1 & \leq N \cdot (s_i - \frac{d_i}{K}) + r \qquad (i>1)
	\tag{\ref{ie6}}\\
e_1 & \leq N \cdot (K d_i - s_i) + r \qquad (i>1)
	\tag{\ref{ie7}}
\end{align}
In the above system, only (\ref{ie51}) is different from (\ref{ie5}) with $m>1$. Inequalities (\ref{ie11}) -- (\ref{ie41}) are special cases of (\ref{ie1}) -- (\ref{ie4}).
Observe that $\lhs{\ref{ie51}} < \lhs{\ref{ie41}} \leq \rhs{\ref{ie41}} < \rhs{\ref{ie51}}$, because $d_1<r$. It follows that Subclaims \ref{wehqroiqhr} -- \ref{q2rugiyi} remain true also in case $m=1$. That is,
$\lhs{\ref{ie11}} \leq \rhs{\ref{ie51}}$, $\lhs{\ref{ie51}} \leq \rhs{\ref{ie11}}$, $\lhs{\ref{ie31}} \leq \rhs{\ref{ie51}}$, $\lhs{\ref{ie51}} \leq \rhs{\ref{ie31}}$, $\lhs{\ref{ie41}} \leq \rhs{\ref{ie51}}$ and $\lhs{\ref{ie51}} \leq \rhs{\ref{ie41}}$. It remains to check the inequalities involving $\rhs{\ref{ie6}}$ and $\rhs{\ref{ie7}}$.

\begin{subclaim}\label{qwriopjqwpr}
$\lhs{\ref{ie11}} \leq \rhs{\ref{ie6}}$, i.e. $e_{1,j} - K d_j  \leq N \cdot (s_i - \frac{d_i}{K}) + r$ for $i,j > 1$.
\end{subclaim}

\begin{pf}
Knowing that $s_i \geq d_i - s_{i,i}$ and $N \cdot s_{i,i} \leq r-e_{i,1}$ (see (\ref{bidobroc})), observe that
\begin{equation}
N \cdot \left(s_i - \frac{d_i}{K}\right) + r \geq N\left(1-\frac 1K\right)\cdot d_i + e_{i,1}.
	\tag{1}\label{o3ut1}
\end{equation}
Thus, using the fact that $e_{1,j}-e_{i,1}\leq e_{i,j}\leq K d_{i,j}$, inequality $\lhs{\ref{ie11}} \leq \rhs{\ref{ie6}}$ is implied by
\begin{equation}
K d_{i,j} - K d_j \leq N\left(1-\frac 1K\right)\cdot d_i.
	\tag{2}\label{o3ut2}
\end{equation}
Further, $K d_{i,j} - K d_j \leq K d_i$, therefore (\ref{o3ut2}) is implied by.
\begin{equation}
K\leq N\left(1-\frac1K\right).
	\tag{3}\label{o3ut3}
\end{equation}
Finally, (\ref{o3ut3}) is satisfied if and only if $\frac{N - \sqrt{N^2 - 4N}}{2} \leq K \leq
\frac{N + \sqrt{N^2 - 4N}}{2}$, which is one of our assumptions.
\end{pf}

\begin{subclaim}
$\lhs{\ref{ie31}} \leq \rhs{\ref{ie6}}$, i.e. $\frac{1}{K} \cdot d_1 \leq N \cdot (s_i - \frac{d_i}{K}) + r$.
\end{subclaim}

\begin{pf}
Using inequality (\ref{o3ut1}) from the proof of Subclaim \ref{qwriopjqwpr} together with $1-\frac1K\geq \frac1N$ (Proposition \ref{dwadwa}), we see that our inequality is implied by
$\frac 1K \cdot d_1 \leq d_i + e_{i,1}$.
This, by the fact that $f$ is $K$-bilipschitz, is implied by
$\frac 1K d_1 \leq d_i + \frac 1K \cdot d_{i,1}$.
The last inequality is true, because $d_1 - d_{i,1} \leq d_i$ and $\frac 1K<1$.
\end{pf}

\begin{subclaim}
$\lhs{\ref{ie41}} \leq \rhs{\ref{ie6}}$, i.e. $(1 + \frac{1}{N}) \cdot d_1 - \frac{1}{N} \cdot r \leq N \cdot (s_i - \frac{d_i}{K}) + r$.
\end{subclaim}

\begin{pf}
Knowing inequalities $s_i \geq d_i - s_{i,i}$, $-N s_{i,i} \geq d_{i,1}-r$ (see (\ref{bidobroc})) and $N(1-\frac1K)\geq1$ (see Proposition \ref{dwadwa}), we obtain
$$
N \cdot \left(s_i - \frac{d_i}{K}\right) + r \geq N\left(1-\frac 1K\right)\cdot d_i + d_{i,1} 
\geq d_i + d_{i,1} \geq d_1.
$$
Finally, $d_1 \geq (1 + \frac1N)\cdot d_1 - \frac 1N\cdot r$, because $d_1<r$.
\end{pf}

\begin{subclaim}\label{pouqrpr}
$\lhs{\ref{ie11}} \leq \rhs{\ref{ie7}}$, i.e. $e_{1,j} - K d_j \leq N \cdot (K d_i - s_i) + r$, \ $j,i > 1$.
\end{subclaim}

\begin{pf}
Using (\ref{bidobroc}), some triangle inequalities and $e_{i,j}\leq K d_{i,j}$, we get
\begin{align*}
N \cdot (K d_i - s_i) + r &=
(NK - N) \cdot d_i - N \cdot (s_i - d_i) + r
\geq (NK - N) \cdot d_i - N \cdot s_{i,i} + r\\
&\geq (NK - N) \cdot d_i - (r - e_{i,1}) + r
= (NK - N) \cdot d_i + e_{i,1}\\
&\geq (NK - N) \cdot d_i + e_{1,j} - e_{i,j} \geq (NK - N) \cdot d_i + e_{1,j} - K d_{i,j}\\
&\geq N(K - 1) \cdot d_i + e_{1,j} - K d_j - K d_i.
\end{align*}
It remains to check that $N(K - 1) \cdot d_i - K d_i \geq0$. This follows from the fact that $N\geq 1 + \frac{1}{K-1}$ (see Proposition \ref{dwadwa}).
\end{pf}

\begin{subclaim}
$\lhs{\ref{ie31}} \leq \rhs{\ref{ie7}}$, i.e. $\frac{1}{K} \cdot d_1 \leq N \cdot (K d_i - s_i) + r$.
\end{subclaim}

\begin{pf}
Repeating the beginning of the proof of Subclaim \ref{pouqrpr} and using the fact that $d_{i,1} \leq K e_{i,1}$, we get
$$
N \cdot (K d_i - s_i) + r \geq (NK - N) \cdot d_i + e_{i,1}
\geq (NK - N) \cdot d_i +
\frac{1}{K} \cdot d_{i,1}.
$$
By Proposition \ref{dwadwa}, $NK-N\geq 1\geq \frac 1K$, therefore finally
$\rhs{\ref{ie7}}\geq \frac1K \cdot d_i + \frac1K \cdot d_{i,1}\geq \frac 1K\cdot d_i = \lhs{\ref{ie31}}$.
\end{pf}

\begin{subclaim}
$\lhs{\ref{ie41}} \leq \rhs{\ref{ie7}}$, i.e. $(1 + \frac{1}{N}) \cdot d_1 - \frac{1}{N} \cdot r \leq 
N \cdot (K d_i - s_i) + r$.
\end{subclaim}

\begin{pf}
Using (\ref{bidobroc}) and some triangle inequalities, we have
\begin{align*}
\rhs{\ref{ie7}} - \lhs{\ref{ie41}} &= NK d_i  + N(d_i-s_i) - N\cdot d_i + (1+\frac1N)(r-d_1)\\
&\geq N(K-1)d_i - N s_{i,i} + (r-d_1) \geq N(K-1)d_i -(r-d_{i,1}) + r - d_1\\
&\geq N(K-1)d_i - d_i = (NK-N-1) d_i.
\end{align*}
Finally, $NK-N-1\geq0$, because $K\geq 1+\frac1N$.
\end{pf}

We have checked all required inequalities, thus showing that the system of inequalities (\ref{ie1}) -- (\ref{ie5}) plus (\ref{ie6}), (\ref{ie7}) has a solution. Thus, the inductive procedure of finding distances $e_m$ can be carried out. This completes the proof of Claim \ref{klejm1}.
\end{pf}

\begin{pf}[Proof of Theorem \ref{t2.1}]
Let $\setof{x_n}{n \in \Nat}$ be a dense subset of $B(x_1,r)$.
We define by induction a sequence of finite functions
$\setof{f_n}{n \in \Nat}$. Let $f_0 = f$.
Suppose that $f_n$ has been defined. Assume that $f_n$ is
$(K,N)$-compliant. By Lemma \ref{l2.4}, there is a $(K,N)$-compliant
function $g_n \sups f_n$ such that $x_n \in \dom(g_n)$.
And there is a $(K,N)$-compliant
function $f_{n + 1} \sups g_n$ such that
$x_n \in \rng(f_{n + 1})$.
Let $h = \bigcup_{n \in \Nat} f_n$.
Since $\dom(h),\rng(h)$ are dense subsets of $B(x_1,r)$ and $h$
is $K$-bilipschitz, there is $g \in \homeo(B(x,r))$
such that $g \sups h$.
It is also obvious that $g$ is $(K,N)$-compliant.
By Proposition~\ref{p2.2},
$g \cup \id{\bbU \sm B(x_1,r)}$ is $K$-bilipschitz.
\end{pf}

\section{Extending uniformly continuous functions}

\index{bicontinuous map}
We say that $\map fXY$ is {\em $(\beta,\al)$-bicontinuous} if $f$ is a $\beta$-continuous bijection such that $f^{-1}$ is $\al$-continuous. In other words, a bijection $\map fXY$ is $(\beta,\al)$-bicontinuous if
$$\al^{-1}(d(x,x')) \loe d(f(x),f(x')) \loe \beta(d(x,x'))$$ holds for every $x,x'\in X$.
Note that this makes sense only if $\al^{-1}\loe \beta$, at least on the range of the metric $d$ of $X$. It turns out that this is not sufficient for the existence of $(\beta,\al)$-bicontinuous extensions.

Recall that every modulus of continuity $\al$ satisfies $$\al(s+t)\loe \al(s)+\al(t)\qquad\text{for every }\,s,t\goe0.$$

\begin{lm}\label{opwtetmsgdpow} Let $\al,\beta\in\mc$. Assume $X\cup\sn p$, $Y$ are finite metric spaces and $\map fXY$ is $(\beta,\al)$-bicon\-ti\-nuous, where $\al,\beta\in\mc$ are such that
\begin{equation}\al^{-1}(s)+\beta(t)\goe \al^{-1}(s+t)\qquad\text{for every }\,s,t\goe0. \tag{*}\end{equation}
Assume $q\notin Y$. Then the formula 
\begin{equation}
d(q,y)=\min\setof{d(f(z),y)+\beta(d(z,p))}{z\in X}.
\tag{**}\end{equation}
defines a metric extension $Y\cup\sn q$ of $Y$ such that $f\cup\sn{\pair pq}$ is $(\beta,\al)$-bicontinuous.
\end{lm}

\begin{pf} In order to justify that (**) defines a metric on $Y\cup\sn q$, we use an argument from \cite{KS}. Define a two-place symmetric function $\phi$ on $Y\cup\sn q$ by setting $\phi(y_0,y_1)=d(y_0,y_1)$ for $y_0,y_1\in Y$ and $\phi(q,y)=\beta(d(p,f^{-1}(y)))$. Then the formula
$$\ovr\phi(y,z)=\min\bigsetof{\sum_{i=0}^{k-1}\phi(y_i,y_{i+1})}{y=y_0,y_1,\dots,y_k=z\text{ and }k\in\omega}$$
clearly defines a metric, called the {\em shortest path metric} measured by $\phi$. It is straight to see that $\ovr\phi(q,y)=d(f(z),y)+\beta(d(z,p))$, i.e. the shortest path from $q$ to $y$ is of the form $(q,f(z),y)$ for some $z\in X$. Note that
$$\phi(f(x_0),q)+\phi(q,f(x_1))\goe\beta(d(x_0,p)+d(p,x_1))\goe\beta(d(x_0,x_1))\goe d(f(x_0),f(x_1)),$$ 
therefore $\ovr\phi(y_0,y_1)=d(y_0,y_1)$ for $y_0,y_1\in Y$. This shows that (**) indeed defines a metric on $Y\cup\sn q$ which extends the metric of $Y$.

Now observe that $d(q,f(x))\loe d(f(x),f(x))+\beta(d(p,x))=\beta(d(p,x))$, i.e. $f\cup\sn{\pair pq}$ is $\beta$-continuous.
Fix $x\in X$ and fix $z\in X$ such that $d(q,f(x))=d(f(x),f(z))+\beta(d(p,z))$. Knowing that $d(f(x),f(z))\goe \al^{-1}(d(x,z))$ and using (*), we get
$$d(q,f(x))\goe \al^{-1}(d(x,z))+\beta(d(p,z))\goe \al^{-1}(d(x,z)+d(p,z))\goe \al^{-1}(d(p,x)).$$
Thus $f^{-1}$ is $\al$-continuous.
\end{pf}

It is easy to see that the above lemma holds for arbitrary (not necessarily finite) metric spaces. The main change in the proof is replacing ``$\min$" by ``$\sup$" in the definition of $\ovr\phi$. For applications, we need the finite version only.

\begin{uwgi}
Let us see that the assumption (*) on $\al,\beta$ is necessary for the existence of extensions. For fix $\al,\beta\in \mc$ such that $\al^{-1}\loe \beta$ and let $X=\{x_0,x_1\}$, $Y=\{y_0,y_1\}$, where $d(x_0,x_1)=s>0$ and $d(y_0,y_1)=\al^{-1}(s)$.
Define $\map fXY$ by $f(x_i)=y_i$. Then $f$ is $(\beta,\al)$-bicontinuous, because $\al^{-1}\loe\beta$. Now fix $p\notin X$ and define $d(p,x_0)=t$ and $d(p,x_1)=s+t$. This defines a metric on $X\cup\sn p$. Suppose $Y\cup\sn q$ is a metric extension of $Y$ such that $f\cup\sn{\pair pq}$ is $(\beta,\al)$-bicontinuous.
Then $d(q,y_0)\loe\beta(t)$ and $d(q,y_1)\goe\al^{-1}(s+t)$. Thus
$$\al^{-1}(s+t)\loe d(q,y_1)\loe d(q,y_0)+d(y_0,y_1)\loe \beta(t)+\al^{-1}(s),$$
which shows that (*) holds.
\end{uwgi}

We shall say that $\beta,\al\in \mc$ are {\em compatible} if 
$$\al^{-1}(s)+\beta(t)\goe \al^{-1}(s+t)\quad\text{ and }\quad \beta^{-1}(s)+\al(t)\goe \beta^{-1}(s+t)$$
holds for every $s,t\goe0$. The above lemma clearly implies extension property for finite $(\beta,\al)$-bicontinuous maps, where $\al,\beta\in\mc$ are compatible. We state this result below.

\begin{wn}\label{c3.3}
Assume $\al,\beta\in \mc$ are compatible moduli of continuity. Then every finite $(\beta,\al)$-bicontinuous bijection between subsets of the Urysohn space $\U$ 
can be extended to a $(\beta,\al)$-bicontinuous homeomorphism of $\U$.
\end{wn}

We now prove a more general version, which involves totally bounded sets. 
The version for isometries was proved by Huhunai\v svili \cite{Hu} in 1955.

\begin{tw}\label{t3.4}
Assume $X,Y\subs\U$ are totally bounded sets and $\map fXY$ is a $(\beta,\al)$-bicontinuous map, where $\al,\beta$ are compatible moduli of continuity. Then there is a $(\beta,\al)$-bicontinuous map $\map F\U\U$ which extends $f$. 
\end{tw}

\begin{pf}
We may assume that both $X,Y$ are closed, since $f$ has a unique continuous extension onto the closure of $X$.
Since the assumptions are symmetric, it suffices to show that $f$ can be extended by adding one point to its domain. Then, by the separability of $\U$, a standard back-and-forth argument will complete the proof.

Fix $p\in\U\sm X$. Fix a sequence $\eps_0>\eps_1>\dots>0$ such that $\beta(\eps_n)\loe 2^{-n}$ for every $n\in\omega$. For each $n\in\omega$ choose an $\eps_n$-net $D_n\subs X$ in such a way that $D_{n+1}\sups D_n$ for every $n\in\omega$.
Let $E_n=\img f{D_n}$ and $f_n=f\rest D_n$.  
We construct inductively a sequence $\sett{q_n}{n\in\omega}\subs\U$ such that $f_n\cup\sn{\pair p{q_n}}$ is $(\beta,\al)$-bicontinuous and $d(q_n,q_{n+1})<2^{-n+1}$.

Start with any $q_0$ obtained by applying Lemma \ref{opwtetmsgdpow} and by the ultrahomogeneity of $\U$. Note that in the construction we will need to use formula (**) given by this lemma.

Now suppose $q_n$ has been already constructed. Apply Lemma \ref{opwtetmsgdpow} to get a metric extension $Z=E_{n+1}\cup\sn q$ of $E_{n+1}$ such that $f_{n+1}\cup\sn{\pair pq}$ is $(\beta,\al)$-bicontinuous and the metric on $Z$ is given by (**). Now $Z$ and $E_{n+1}\cup\sn{q_n}$ are two compatible metric spaces whose intersection is $E_{n+1}$. We can amalgamate them in the minimal way, i.e. setting
$$d(q,q_n)=\max\setof{|d(q,y)-d(y,q_n)|}{y\in E_{n+1}}.$$
Let us now estimate $d(q,q_n)$ from above. Fix $x\in D_{n+1}$ such that 
$$d(q,q_n)=|d(q,f(x))-d(f(x),q_n)|.$$
Applying (**) we see immediately that $d(q,f(x))\loe d(q_n,f(x))$, because $D_n\subs D_{n+1}$. 
Now find $u\in D_{n+1}$ with $$d(q,f(x))=d(f(x),f(u))+\beta(d(u,p))$$ and
find $z\in D_n$ such that $d(u,z)<\eps_n$. Then
$$d(q_n,f(x))\loe d(f(x),f(z))+\beta(d(z,p))$$
and hence
\begin{align*}
d(q_n,f(x))-d(q,f(x))&\loe d(f(x),f(z))+\beta(d(z,p)) - d(f(x),f(u))-\beta(d(u,p))\\
&\loe d(f(z),f(u)) + \Bigl(\beta(d(z,p))-\beta(d(u,p))\Bigr)\\
&\loe \beta(d(z,u)) + \Bigl(\beta(d(z,u)+d(u,p))-\beta(d(u,p))\Bigr)\\
&< \beta(\eps_n) + \beta(d(z,u))+\beta(d(u,p))-\beta(d(u,p))\\
&< 2\beta(\eps_n)<2^{-n+1}.
\end{align*}
Thus we have proved that $d(q,q_n)=|d(q_n,f(x))-d(q,f(x))|<2^{-n+1}$.

Finally, we find $q_{n+1}\in \U$ which realizes our amalgamation, i.e. $d(q_{n+1},q_n)=d(q,q_n)$. This finishes the description of the inductive construction.

Clearly, $\sett{q_n}{n\in\omega}$ is a Cauchy sequence in $\U$. Let $q=\lim_{n\to\infty}q_n$. We claim that $f\cup\sn{\pair pq}$ is $(\beta,\al)$-bicontinuous. Indeed, given $x\in D_k$ and $n\goe k$ we have
$$d(q,f(x))\loe d(q,q_n)+d(q_n,f(x))<2^{-n+1}+\beta(d(p,x))$$
and
$$d(q,f(x))\goe d(q_n,f(x))-d(q_n,q)> \al^{-1}(d(p,x))-2^{-n+1}.$$
Passing to the limit, we get
$$\al^{-1}(d(p,x))\loe d(q,f(x))\loe\beta(d(p,x)).$$
Since $\bigcup_{n\in\omega}D_n$ is dense in $X$, the above inequalities hold for every $x\in X$. This completes the proof.
\end{pf}

It turns out that the Bilipschitz Extension Theorem can be generalized to the case of totally bounded subsets of a ball, using ideas from the above proof and elaborating arguments from the proof of Lemma \ref{l2.4}. The precise statement looks as follows.

% TOTALLY BOUNDED VERSION

\begin{tw}\label{t2.91}
Let $B(x,r)$ be an open ball in $\bbU$,
let $N \geq 4$, $\frac{N - \sqrt{N^2 - 4N}}{2} < K < \frac{N + \sqrt{N^2 - 4N}}{2}$ and
let $A \subs B(x,r)$ be a totally bounded set. Assume further that $\map{f}{A}{B(x,r)}$ is $K$-bilipschitz and $N$-bigood and $x\in A$ is such that $f(x)=x$.

Then there exists a bijection
$\map{g}{B(x,r)}{B(x,r)}$ such that
\begin{itemize}
\item[\num{1}]
$g \sups f$,
\item[\num{2}]
$g \cup \id{\bbU \sm B(x,r)}$
is $K$-bilipschitz,
\item[\num{4}]
$g$ is $N$-bigood.
\end{itemize}
\end{tw}

\section{A metric on the bilipschitz group}

In order to prove our main result, we need to know that the group of bilipschitz auto-homeomorphisms of a metric space can be endowed with a suitable metrizable topology, compatible with the group structure. This is the contents of the current section. The results are rather standard, however we were unable to find any bibliographic references, therefore we give all the details.

Let $\pair{X}{d}$ be a metric space.
Recall that $\lip(\pair{X}{d})$ denotes the group of the bilischitz
auto-homeomorphisms of $\pair Xd$. For $g \in \lip(X)$ let
$$\lipc(g) = \min\setof{K}{g \mbox{ is $K$-bilipschitz}}.$$
We define the following semimetrics on $\lip(X)$.
$$d_L(f,g) = \log(\lipc(f^{-1}\cmp g)).$$
Let $x_0 \in X$. For $n \in \Nat$ define
$$d_n(f,g) = \sup\setof{d(f(x),g(x))}{x \in B(x_0,n)}.$$
Further, define
$$d_S(f,g) = \sum_{n \in \Nat} \frac{d_n(f,g)}{2^n}$$
and
$$\hat d(f,g) = \max(d_L(f,g),d_S(f,g)).$$
Let $\tau_L$ be the topology of the semimetric $d_L$
and $\tau_n$ be the topology of the semimetric $d_n$.
Also, $\tau$ will denote the topology of $d$. Finally, $\hat \tau$ will be the topology on $\lip(X)$ induced by $\hat d$.

\begin{prop}
Let $\pair Xd$ be as above.
\begin{enumerate}
	\item[(a)] For every $f,g \in \lip(X)$, $d_S(f,g) < \infty$.
	\item[(b)] $\hat d$ is a metric on $\lip(X)$.
	\item[(c)] The topology $\hat\tau$ is generated by $\tau_L \cup \bigcup_{n \in \Nat} \tau_n$.
\end{enumerate}
\end{prop}

\begin{pf}
Let $G = \lip(X)$.

(a) Let $f,g \in G$ be $K$-Lipschitz and fix $n \in \Nat$ and $x \in B(x_0,n)$.
Denote $a = d(f(x_0),g(x_0))$. Then
$$d(f(x),g(x)) \leq d(f(x),f(x_0)) + d(f(x_0),g(x_0)) + d(g(x_0),g(x)) \leq
= 2Kn +a.$$
So $d_S(f,g) \leq \sum_{n \in \Nat} \frac{2Kn + a}{2^n}$ and the series is convergent.

Part (b) is trivial.

(c) It is obvious that for every $n \in \Nat$,
$\tau_n \subs \hat\tau$ and $\tau_L \subs \hat\tau$.
Let $f \in G$ and $r > 0$.
Denote $B = B^{d_S}(f,r)$.
Let $B_L = B^{d_L}(f,\log(2))$ and let $B_1 = B^{d_1}(f,\frac{r}{2})$.
Denote $K = \lipc(f)$. Let $g \in B_L \cap B_1$.
Then
$\lipc(g) \leq \lipc(f) \cdot \lipc(f^{-1}\cmp g) < 2K$.

Suppose that $x \in B(x_0,n)$.
Then
$$d(f(x),g(x)) \leq 
d(f(x),f(x_0)) + d(f(x_0),g(x_0)) + d(g(x_0),g(x)) <
Kn + \frac{r}{2} + 2Kn = 3Kn + \frac{r}{2}.$$
That is,
\begin{equation}
d(f(x),g(x))  < 3Kn + \frac{r}{2}.
\tag{1}\label{equ1}
\end{equation}
There is $n_0$ such that
\begin{equation}
\sum_{n > n_0}\frac{3Kn + \frac{r}{2}}{2^n} < \frac{r}{4}.
\tag{2}\label{equ2}
\end{equation}
Let $s > 0$ be such that
$\sum_{n = 2}^{n_0}\frac{s}{2^n} < \frac{r}{4}$.
Let
$C = B_L \cap B_1 \cap \bigcap_{n = 2}^{n_0} B^{d_n}(f,s)$.
We show that $C \subs B$. Let $g \in C$.
Since $g \in B_1$,
\begin{equation}
d_1(f,g) < \frac{r}{2}.
\tag{3}\label{equ3}
\end{equation}
Since $g \in \bigcap_{n = 2}^{n_0} B^{d_n}(f,s)$,
\begin{equation}
\sum_{n = 2}^{n_0}\frac{d_n(f,g)}{2^n} < \frac{r}{4}.
\tag{4}\label{equ4}
\end{equation}
Let $n > n_0$ and $x \in B(x_0,n)$.
By (\ref{equ1}),
$d(f(x),g(x))  < 3Kn + \frac{r}{2}$.
Hence $d_n(f,g) \leq 3Kn + \frac{r}{2}$.
By (\ref{equ2}), we get

\begin{equation}
\sum_{n > n_0}\frac{d_n(f,g)}{2^n} < \frac{r}{4}.
\tag{5}\label{equ5}
\end{equation}
It follows from (\ref{equ3}) -- (\ref{equ5}) that $d_S(f,g) < r$.
That is, $g \in B$.

We have shown that every ball $B$ in the metric $d_S$ with center
at $f$ contains a finite intersection of balls with center $f$ in the
semimetrics $d_L,d_1,\dots,d_n,\ldots$.
It follows that $B$ is open in the topology generated
by $\tau_L \cup \bigcup_{n \in \Nat} \tau_n$.
\end{pf}

\begin{tw}\label{t4.2}
Let $\pair Xd$ be a metric space. Then:
\begin{enumerate}
	\item[(a)] $\pair{\lip(X)}{\hat\tau}$ is a topological group.
	\item[(b)] The action of $\pair{\lip(X)}{\hat\tau}$ on $\pair{X}{\tau}$ is continuous.
	\item[(c)] If $\pair Xd$ is complete, then so is $\pair{\lip(X)}{\hat d}$.
\end{enumerate}
\end{tw}

\begin{pf} Let $G = \lip(X)$.

(a) Let $n \in \Nat$. Let $f,g \in G$
and let $\varepsilon > 0$.
Let $K = \lipc(f)$.
Let $m$ be such that $g(B(x_0,n)) \subs B(x_0,m)$
and let $\ell > 2nK$.
Take any
$$f_1 \in B^{d_m}(f,\frac{\varepsilon}{2}) \cap B^{d_L}(f,\log(2)) \oraz g_1 \in B^{d_n}(g,\frac{\varepsilon}{4K}).$$
So $\lipc(f_1) < 2K$.
Let $x \in B(x_0,n)$. So $g(x) \in B(x_0,m)$. We have
\begin{align*}
d(f_1g_1(x),fg(x)) &\leq
d(f_1g_1(x),f_1g(x)) + d(f_1g(x),fg(x))\\ &<
2K \cdot d(g_1(x),g(x)) + \frac{\varepsilon}{2}<
2K \frac{\varepsilon}{4K} + \frac{\varepsilon}{2} = \varepsilon.
\end{align*}
This shows that the inverse image of a $d_S$-open subset of $G$
under multiplication is $\hat d$-open in $G \times G$.

We now show that the inverse image of a $d_L$-open subset of $G$ under
multiplication is $\hat d$-open in $G \times G$.
Let $f,g \in G$ and $M > 1$.
We shall show that there are $d_L$-neighborhoods $U,V$ of $f$ and $g$
respectively such that if $f_1 \in U$ and $g_1 \in V$,
then $\lipc((f\cmp g)^{-1} f_1\cmp g_1) < M$.
Let $K = \lipc(g)$ and let
$$U = \setof{h \in G}{\lipc(f^{-1}\cmp h) < \frac{\sqrt{M}}{K^2}},\quad V = \setof{h \in G}{\lipc(g^{-1}\cmp h) < \sqrt{M}}.$$
Fix $f_1 \in U$ and $g_1 \in V$. Then
\begin{align*}
\lipc((f\cmp g)^{-1} f_1\cmp g_1) &=
\lipc(g^{-1}\cmp f^{-1} f_1\cmp g_1)\\ &=
\lipc((g^{-1}\cmp f^{-1}\cmp f_1\cmp g)\cmp (g^{-1}\cmp g_1))\\ &<
K \frac{\sqrt{M}}{K^2} K \cdot \sqrt{M} = M.
\end{align*}
So $U \cdot V \subs B^{d_L}(f\cmp g, \log(M))$.

The function $f \mapsto f^{-1}$ is an isometry of $d_L$.
It thus remains to show that for every $n \in \Nat$,
$f \in G$ and $V \in \tau^{d_n}$ such that $f^{-1} \in V$,
there exists a neighborhood $U$ of $f$ satisfying $U^{-1} \subs V$.

Suppose that $V = B^{d_n}(f^{-1},\varepsilon)$.
Let $K = \lipc(f)$.
Let $m$ be such that $f^{-1}(B(x_0,n)) \subs B(x_0,m)$.
Let $\ell > m + 2 \varepsilon$.
We shall show that
$$U = B^{d_{\ell}}(f,\frac{\varepsilon}{K}) \cap B^{d_L}(f,\log(2)).$$
Let
$f_1 \in
B^{d_{\ell}}(f,\frac{\varepsilon}{K}) \cap B^{d_L}(f,\log(2))$.
So $\lipc(f_1) < 2K$.
Let $x \in B(x_0,n)$. Denote $y = f^{-1}(x)$ and $z = f_1^{-1}(x)$.
We prove that $z \in B(x_0,\ell)$.
Note that $y \in B(x_0,m)$, so $d(f(y),f_1(y)) < \frac{\varepsilon}{K}$.
That is,
$$d(f_1(z),f_1(y)) < \frac{\varepsilon}{K}.$$
Hence $d(y,z) < 2 \varepsilon$ and therefore
$d(x_0,z) \leq d(x_0,y) + d(y,z) < m + 2 \varepsilon < \ell$.
Thus $z \in B(x_0,\ell)$.

It follows that $d(f(z),f_1(z)) < \frac{\varepsilon}{K}$.
That is, $d(f(z),f(y)) < \frac{\varepsilon}{K}$.
Since $f$ is $K$-bilipschitz,
$d(y,z) < \varepsilon$.
Finally, $f_1 \in B^{d_n}(f,\varepsilon)$.

(b) Suppose that $f(x) = y$ and let $\varepsilon > 0$.

Let $\delta$ be such that if $d(x_1,x) < \delta$,
then 
$d(f(x_1),f(x)) < \frac{\varepsilon}{2}$.
Let $n$ be such that $B(x,\delta) \subs B(x_0,n)$.
Let
$x_1 \in B(x,\delta)$
and
$f_1 \in B^{d_n}(f,\frac{\varepsilon}{2})$.
Then
$d(f_1(x_1),f(x)) \leq
d(f_1(x_1),f(x_1)) + d(f(x_1),f(x)) \leq
\frac{\varepsilon}{2} + \frac{\varepsilon}{2} = \varepsilon$.
So the action of $G$ on $X$ is continuous.

(c) Let $\setof{f_n}{n \in \Nat}$ be a Cauchy sequence.
Then for every $m \in \Nat$, $\sett{f_n}{n \in \Nat}$
is a Cauchy sequence with respect to $d_m$.
So there is a continuous function $g_m$ such that
$g_m$ is the uniform limit of
$\setof{f_n \rest B(x_0,m)}{n \in \Nat}$.
Let $g = \bigcup_{m \in \Nat} g_m$.
Since $\setof{f_n}{n \in \Nat}$ is a Cauchy sequence, it is bounded.
That is, there is $r$ such that
$\setof{f_n}{n \in \Nat} \subs B^{d_L}(\id{},r)$.
Hence there is $K$ such that for every $n \in \Nat$, $f_n$ is $K$-bilipschitz.
$g$ is the pointwise limit of $\setof{f_n}{n \in \Nat}$. So $g$ is bilipschitz.

For every $m \in \Nat$, $g$ is a $\tau_m$ limit of $\setof{f_n}{n \in \Nat}$.
It remains to show that $g$ is a $\tau_L$ limit of
$\sett{f_n}{n \in \Nat}$.
Let $k > 1$. There is $n_0$ such that for every $n,m \geq n_0$
$\lipc(f_n^{-1}\cmp f_m) < k$.
Let $n \geq n_0$. $f_n^{-1}\cmp g$ is the pointwise limit of
$\sett{f_n^{-1}\cmp f_m}{m \in \Nat}$.
Since for every $m$, $f_n^{-1}\cmp f_m$ is $k$-bilpschitz,
$f_n^{-1}\cmp g$ is $k$-bilpschitz.
It follows that for every $n \geq n_0$,
$\lipc(f_n^{-1}\cmp g) \leq k$.
So $g$ is a $\tau_L$ limit of the sequence $\sett{f_n}{n \in \Nat}$.
\end{pf}

\section{Obtaining a homeomorphism from a group isomorphism}
\label{s5}

Suppose that $X$ and $Y$ are open subsets of $\bbU$,
and
$\iso
{\varphi}{\homeolc{\Gamma}(X)}{\homeolc{\Delta}(Y)}$.
We shall show that there is $\iso{\tau}{X}{Y}$
such that $\varphi(g) = \tau \cmp g \cmp \tau^{-1}$
for every $g \in \homeolc{\Gamma}(X)$.
The proof relies on a theorem from \cite{FR}.
In order to state it, we introduce some new notions.

\index{space-group pair}\index{faithful class}
Let $X$ be a topological space and $G$ be a subgroup of
the group $\homeo(X)$ of all auto-homeo\-mor\-phisms of $X$.
The pair $\pair{X}{G}$ is then called a
{\it space-group pair}.
Let $K$ be a class of space-group pairs.
$K$ is called a {\it faithful class} if for every
$\pair{X}{G}, \pair{Y}{H} \in K$
and an isomorphism $\varphi$ between the groups
$G$ and $H$ there is a homeomorphism $\tau$ between $X$ and $Y$
such that $\varphi(g) = \tau \cmp g \cmp \tau^{-1}$
for every $g \in G$.
Let $\pair{X}{G}$ be a space-group pair and
$S \subs X$ be open. $S$ is {\it strongly flexible},
if for every infinite $A \subs S$ without
accumulation points in $X$,
there is a nonempty open set $V \subs X$ such that
for every nonempty open set $W \subs V$ there is $g \in G$
such that the sets
$\setof{a \in A}{g(a) \in W}$ and
$\setof{a \in A}{\mbox{for some neighborhood } U \mbox{ of } a,
g \rest U  = \id{}} $
are infinite.
\index{strongly flexible set}
\index{set!-- strongly flexible}

\begin{tw}[\cite{FR}, Theorem B]\label{t5.1}
Let $K_F$ be the class of all space-group pairs $\pair{X}{G}$
such that
\begin{itemize}
	\item[{(1)}] $X$ is regular, first countable and has no isolated points.
	\item[{(2)}] For every $x \in X$ and an open neighborhood $U$ of $x$ the set
$$\setof{g(x)}{g \in G \mbox{ and } g \rest (X\sm U) = \id{} }$$ is somewhere dense.
	\item[{(3)}] The family of strongly flexible sets is a cover of $X$.
\end{itemize}
Then $K_F$ is faithful.
\end{tw}

We wish to show that if $X$ is an open subset of $\bbU$
then
$\pair{X}{\homeolc{\Gamma}(X)} \in K_F$.
Note that if $\pair{X}{G} \in K_F$ and $G \leq H \leq \homeo(X)$,
then $\pair{X}{H} \in K_F$.
So it suffices to show that if $X$ is an open subset of
$\bbU$, then $\pair{X}{\homeolc{\sLip}(X)} \in K_F$.

Clause (1) in the definition of $K_F$ certainly holds for
open subsets of $\bbU$,
and Clause (2) follows trivially from Theorem~\ref{t2.1}.
So it remains to show that open subsets of $\bbU$
have a cover consisting of strongly flexible sets.
The proof of this fact is the contents of this section.

%\separator

Suppose that $\pair{X}{G}$ is a space-group pair
and $A \subs X$ is infinite.
\index{dissectable set}\index{set!-- dissectable}
We say that $A$ is {\em dissectable} with respect to $\pair{X}{G}$,
if there is a nonempty open set $V \subs X$ such that
for every nonempty open set $W \subs V$ there is $g \in G$
such that the sets
$$\setof{a \in A}{g(a) \in W} \oraz 
\setof{a \in A}{g \rest U  = \id{}\mbox{ for some neighborhood } U \mbox{ of }a}$$
are infinite.

Let $X$ be a topological space and
$A \subs X$. We say that $A$ is {\it completely discrete},
if $A$ has no accumulation points.
\index{set!-- completely discrete}

Suppose that $X$ is a metric space, $A \subs X$
and $r > 0$. We say that $A$ is {\em $r$-spaced}, if $d(a,b) > r$
for every distinct $a,b \in A$.
We say that $A$ is {\em spaced}, if for some $r > 0$,
$A$ is $r$-spaced.
\index{spaced set}\index{set!-- spaced}

It follows immediately from the definition, that
if $A$ is dissectable and $B \sups A$,
then $B$ is dissectable.
Since $\bbU$ is a complete metric space,
every  completely discrete set contains a spaced subset.
For spaces which are open subsets of $\bbU$ we shall prove that
every spaced set contained in a small ball is dissectable.

Suppose that $a < b$ are real numbers,
and $\map{h}{[a,b]}{\bbU}$ is an isometry into $\bbU$.
Then $L := h([a,b])$ is called a
{\em line segment} in $\bbU$,
and $h(a),h(b)$ are the {\em endpoints} of $L$.

Suppose that $X$ is a metric space and for every $n \in \Nat$,
$\map{f_n}{[a,b]}{X}$.
We say that $\setof{f_n(t)}{n \in \Nat}$ is {\it equicontinuous},
if for every $\eps > 0$ there is $\delta > 0$
such that for every $t,s \in [a,b]$ and $n \in \Nat$:
if $\abs{t - s} < \delta$, then $d(f_n(s),f_n(t)) < \eps$.
Let $r > 0$ and $\calA$ be a set of subsets of $X$.
We say that $\calA$ is {\it $r$-spaced} if for every distinct
$A,B \in \calA$, $d(A,B) > r$.
We say that $\calA$ is {\it spaced} if for some $r > 0$,
$\calA$ is $r$-spaced.
We say that $\calA$ is {\it almost $r$-spaced} if for some finite
$\calA_0 \subs \calA$, $\calA\sm \calA_0$ is $r$-spaced.
We say that $\calA$ is {\it almost spaced} if for some $r > 0$,
$\calA$ is almost $r$-spaced.
Let $A \almostcontained B$ mean that $A\sm B$ is finite.

Given a set $B\subs X$ and a family $\Ef$ of self-maps of a space $X$, we shall denote by $\stablep\Ef B$ the set of all $f\in \Ef$ such that $f\rest(X\sm B)=\id{}$.

\begin{prop}\label{p5.6}
(a) Let $x,u,v \in \bbU$ and $r > 0$ be such that
$u,v \in B(x,\frac{r}{15})$.
Then there is $g \in \stablep {\homeo(\bbU)}{B(x,r)}$
such that $g(u) = v$ and $g$ is $2$-bilipschitz.

(b) There is an increasing function
$\map{\hat K}{(0,\infty)}{\Nat}$
such that the following holds.
If $L$ is a line segment in $\bbU$ with endpoints
$x$ and $y$ and $r > 0$, then there is $g \in \stablep {\homeo(\bbU)}{B(L,r)}$
such that $g(x) = y$ and $g$ is
$\hat K(\frac1r{\length(L)})$-bilipschitz.
\end{prop}

\begin{pf}
(a) Set $s = \frac{r}{15}$, $t = 12s$ and let $y \in \bbU$ be such that
$d(x,y) = 3s$.
Then $B(y,12s) \subs B(x,15s) = B(x,r)$.
Let $f = \dn{\pair{y}{y}}{\pair{u}{v}}$.
Then $d(u,y),d(v,y) \in (2s,4s)$ and hence $\lipc(f) < 2$.
Also,
$t - d(u,y),t - d(v,y) > 12s - 4s = 8s$.
So $\frac{d(u,v)}{t - d(u,y)},\frac{d(u,v)}{t - d(v,y)} <
\frac{2s}{8s} = \frac{1}{4}$.
By Theorem~\ref{t2.1}, there is $g \in \stablep{\homeo(\bbU)}{B(y,12t)}$
such that $g \sups f$.
Hence $g(u) = v$ and $\supp(g) \subs B(x,r)$.

(b) Let $\hat K(t) = 2^{[16t] + 1}$.
Choose a sequence of points
$x = x_0,x_1,\dots,x_n = y$ in $L$ such that for $i < n$,
$d(x_{i - 1},x_i) = \frac{r}{16}$
and $d(x_{n - 1},x_n) \leq \frac{r}{16}$.
Then $n \leq [\frac{16 \cdot \length(L)}{r}] + 1$.
Let $s = \frac{r}{15}$. Then for every $i = 1,\dots,n$,
$x_{i - 1} \in B(x_i,s)$ and $B(x_i,15s) \subs B(L,r)$.
By (a), there are $g_1,\dots,g_n \in \stablep{\homeo(\bbU)}{B(L,r)}$
such that $g_i(x_{i - 1}) = x_i$ and $g_i$ is $2$-bilipschitz.
Let $g = g_n \cmp\dots\cmp g_1$.
Then $g(x) = y$, $g \in \stablep{\homeo(\bbU)}{B(L,r)}$
and $g$ is $2^n$-bilipschitz. That is,
$g$ is $\hat K(\frac{\length(L)}{r})$-bilipschitz.
\end{pf}

\begin{prop}\label{p5.8}
Let $X$ be a metric space and $\setof{\gamma_n}{n \in \Nat}$
be a sequence of arcs in $X$ such that
$\map{\gamma_n}{[0,1]}{X}$.
Suppose that
$\setof{\gamma_n(0)}{n \in \Nat}$ is spaced,
and $\setof{\gamma_n(1)}{n \in \Nat}$ is a Cauchy sequence.
Also assume that $\setof{\gamma_n}{n \in \Nat}$ is equicontinuous.

{Then} there are $s \in (0,1]$
and an infinite $\sigma \subs \Nat$
such that $\setof{\gamma_n(s)}{n \in \sigma}$ is a Cauchy sequence,
and for every $t \in [0,s)$,
$\setof{\gamma_n([0,t])}{n \in \sigma}$ is almost spaced.
\end{prop}

\begin{pf}
For every infinite $\eta \subs \Nat$ define
$$
s_{\eta} =
\sup(\setof{s \in [0,1]}
{\setof{\gamma_n([0,s])}{n \in \eta}
\mbox{ is almost spaced}}).
$$
Clearly, if $\eta \almostcontained \zeta$,
then $s_{\eta} \geq s_{\zeta}$.
This implies that there is an infinite $\sigma \subs \Nat$
such that for every infinite $\eta \subs \sigma$,
$s_{\eta} = s_{\sigma}$. Denote $s_{\sigma}$ by $s$.
We show that there is an infinite $\sigma_0 \subs \sigma$
such that $\setof{\gamma_n(s)}{s \in \sigma_0}$ is a Cauchy sequence.
If $s = 1$, then $\setof{\gamma_n(s)}{s \in \sigma}$
is a Cauchy sequence.
Suppose that $s < 1$.
If there is no $\sigma_0$ as required, then there are $r > 0$
and an infinite $\sigma_0 \subs \sigma$
such that $\setof{\gamma_n(s)}{s \in \sigma_0}$ is $r$-spaced.
Let $\delta > 0$ be such that for every $n \in \Nat$ and
$t,u \in [0,1]$: if $\abs{t - u} < \delta$,
then $d(\gamma_n(t),\gamma_n(u)) < \frac{r}{3}$.
We may assume that $s + \delta \leq 1$.
Then $\setof{\gamma_n([0,s + \delta])}{n \in \sigma_0}$
is $\frac{r}{3}$-spaced.
So $s_{\sigma_0} \geq s + \delta$, a contradiction.
Hence, there is an infinite $\sigma_0 \subs \sigma$
such that $\setof{\gamma_n(s)}{n  \in \sigma_0}$ is a Cauchy sequence.
We may thus assume that
$\setof{\gamma_n(s)}{n  \in \sigma}$ is a Cauchy sequence.
By the definition of $s$, for every $t \in [0,s)$,
$\setof{\gamma_n([0,t])}{n \in \sigma}$ is almost spaced.
\end{pf}

\begin{lm}\label{l5.9}
Let $X \subs \bbU$ be open.
Suppose that $r,\eps > 0$,
$B^{\bbU}(x,r + \eps) \subs X$
and $A \subs B^{\bbU}(x,r)$ is an infinite spaced set.
Then $A$ is dissectable with respect to $\stablep{\lip(\bbU)}{X}$.
\end{lm}

\begin{pf}
Suppose that $A = \setof{a_n}{n \in \Nat}$.
For every $n \in \Nat$ let $L_n$ be a line segment connecting
$a_n$ with $x$,
and let $\map{\gamma_n}{[0,1]}{X}$ be the paramet\-rization of $L_n$
such that $d(\gamma_n(t),\gamma_n(0)) = t \cdot d(a_n,x)$.
Then $\setof{\gamma_n}{n \in \Nat}$ is equicontinuous.
Let $s \in [0,1]$ and $\sigma \subs \Nat$ be such that
$\sigma$ is infinite, $\setof{\gamma_n(s)}{n \in \sigma}$
is a Cauchy sequence and for every $t \in [0,s)$,
$\setof{\gamma_n([0,t])}{n \in \sigma}$ is almost spaced.
We may assume that $\sigma = \Nat$.
Let $y = \lim_{n \in \Nat} \gamma_n(s)$.
Clearly, $d(\bigcup_{n \in \Nat} L_n,\bbU \sm X) \geq \eps$.
Hence $d(y,\bbU \sm X) \geq \eps$.
We may assume that $y \not\in \setof{a_n}{n \in \Nat}$.
So $\delta := d(y,\setof{a_n}{n \in \Nat}) > 0$.
Set $s = \frac1{16}{\min(\eps,\delta)}$.
We show that for every nonempty open $W \subs B(y,s)$
there is $g \in \stablep{\lip{\bbU}}{X}$
such that the sets
$$\setof{i \in \Nat}{g(a_i) \in W}
\oraz
\setof{i \in \Nat}{\mbox{there is } V \in \nbd (a_i)
\mbox{ such that } g \rest V = \id{}}$$
are infinite.
We may assume that $W = B(z,q)$.
There is $t < s$ such that $\gamma_n(t) \in B(y,\frac{q}{2})$ for all but finitely many $n$'s.
There are a finite set $\eta \subs \Nat$ and $e > 0$
such that
$\setof{\gamma_n([0,t])}{n \in \Nat \sm \eta}$ is $e$-spaced.
Let $p = \min(\eps,\frac{e}{2})$.
Then for every $m \neq n$ in $\Nat \sm \eta$,
$B(\gamma_m([0,t]),p) \cap B(\gamma_n([0,t]),p) = \emptyset$
and for every $n$, $B(\gamma_n([0,t]),p) \subs X$.
Let $\ell = \sup_{n \in \Nat} \length(L_n)$.
For every $n \in \Nat \sm \eta$ there is
$h_n \in \stablep{\homeo(\bbU)}{B(\gamma_m([0,t]),p)}$ such that
$h_n(a_n) = \gamma_n(t)$ and $h_n$ is
$\hat K(\frac{\ell}{p})$-bilipschitz.
This follows from Propostion \ref{p5.6}(b).
Let
$h = \cmp\; \setof{h_{2n}}{n \in \Nat \mbox{ and } 2n \not\in \eta}$.
It follows trivially from the above
that $h \in \stablep{\homeo(\bbU)}{X}$
and $h$ is $(\hat K(\frac{\ell}{p}))^2$-bilipschitz.
Also for every $n \in (2 \Nat) \sm \eta$,
$h(a_n) \in B(y,\frac{q}{2})$,
and for every $n \in (2 \Nat + 1) \sm \eta$
there is $V \in \nbd (a_n)$
such that $h \rest V = \id{}$.
Note that $d(y,z) < \frac1{16}{\min(\eps,\delta)}$.
By Proposition~\ref{p5.6}(a),
there is
$f \in \stablep{\homeo(\bbU)}{B(y,\frac{15}{16} \cdot \min(\eps,\delta))}$
such that $f(y) = z$ and $f$ is $2$-bilipschitz.
Since $d(y,\setof{a_n}{n \in \Nat}) = \delta$,
for every $n \in \Nat$ there is $V \in \nbd (a_n)$
such that $f \rest V = \id{}$.
It also follows that $f \in \stablep{\homeo(\bbU)}{X}$.
Let $g = f \cmp h$.
Then $g \in \stablep{\homeo(\bbU)}{X}$ and $g$ is bilipschitz.
So
\begin{itemize}
\item[(1)] 
$g \in \stablep{\lip(\bbU)}{X}$.
\end{itemize}
Let $n \in (2 \Nat) \sm \eta$.
Then $h(a_n) \in B(y,\frac{q}{2})$.
Since $f$ is $2$-bilipschitz and $f(y) = z$,
it follows that
$f(h(a_n)) \in B(z,q)$.
That is,
\begin{itemize}
\item[(2)] 
$g(a_n) \in B(z,q)$ for every $n \in (2 \Nat) \sm \eta$.
\end{itemize}
Finally,
\begin{itemize}
\item[(3)] 
For every $n \in (2 \Nat + 1) \sm \eta$
there is $V \in \nbd (a_n)$
such that $g \rest V = \id{}$.
\vspace{-05.7pt}
\end{itemize}
We have shown that $A$ is dissectable.
\end{pf}

\begin{wn}\label{c5.10}
Let $X$ be a nonempty open subset of $\bbU$
and
$\stablep{\lip(\bbU)}{X} \leq G \leq \homeo(X)$.
Then $\pair{X}{G} \in K_F$.
\end{wn}

\begin{pf}
(a) Note that $X$ is a first countable regular space without isolated points.
That is, Clause~1 in the definition of $K_F$ holds.
By Proposition~\ref{p5.6}(a),
for every $x \in X$ and $U \in \nbd^X(x)$
the set
$\setof{g(x)}{g \in \stablep{\lip{\bbU}}{X}}$ is somewhere dense.
So Clause~2 in the definition of $K_F$ holds.

Note that in a complete metric space
every completely discrete infinite set contains an infinite spaced
subset.
It thus follows from Lemma~\ref{l5.9}(a)
that if $B^{\bbU}(x,r + \varepsilon) \subseteq X$,
then $B^{\bbU}(x,r)$ is strongly flexible with respect to
$\stablep{\lip(\bbU)}{X}$,
and thus it is strongly flexible with respect to $G$.
So $X$ has a cover consisting of strongly flexible sets.
That is, Clause~3 holds.
\end{pf}

\section{Local $\Gamma$-bicontinuity of the conjugating homeomorphism}\label{s6}

Let $X \subseteq \bbU$ be open. Define
$$
\lip^{\id{}}(X) =
\setof{g \rest X}{g \in \stablep{\lip(\bbU)}{X}}.
$$
We equip $\lip^{\id{}}(X)$ with the topology it inherits
from $\lip(\bbU)$.
We shall apply Theorem 3.41 from \cite{RY} to the group $\lip^{\id{}}(X)$.
This requires the following definitions.

\begin{df}\label{d6.1}
(a) Let $X$ be a metric space, $G \leq \homeo(X)$
and $\sigma$ be a topology on $G$ such that the action of
$G$ on $X$ is continuous with respect to $\sigma$.
Let $x \in X$.
We say that $\pair{G}{\sigma}$ is {\em affine-like at $x$},
if the following holds.
For every $V \in \nbd(\id{})$ and $W \in \nbd(x)$
there are
$n = n(x,V,W) \in \Nat$,
$U = U(x,V,W) \in \nbd(x)$
such that for every
$x_1,y_1,x_2,y_2 \in U$: if
$d(x_1,y_1) = d(x_2,y_2)$,
then there are $h_1, \dots, h_n \in V$ such that
$h_1 \cmp \dots \cmp h_n(x_1) = x_2$, \,
$h_1 \cmp \dots \cmp h_n(y_1) = y_2$
and $h_i \cmp h_{i + 1} \cmp \dots \cmp h_n(\dn{x_1}{y_1}) \subs W$
for every $1 \leq i \leq n$.
\index{affine-like at $x$}

If $\pair{G}{\sigma}$ is affine-like at every $x \in X$,
then $\pair{G}{\sigma}$
is said to be an {\em affine-like topological group}.
\index{affine-like topological group}

(b) Let $X$ be a metric space and $x \in X$.
We say that $X$ has the {\em discrete path property at $x$} (briefly: {\em $X$ is DPT at $x$}),
if the following holds.
There is $U \in \nbd(x)$ and $K \geq 1$ such that
\begin{itemize}
\item[($*$)] for every $y,z \in U$ and $d \in (0, d(y,z))$
there are $n \in \Nat$ and $u_0,\dots,u_n \in X$
such that $n \leq K \cdot \frac{d(y,z)}{d}$,
$d(y,u_0),\,d(u_n,z) < d$
and \,$d(u_{i - 1},u_i) = d$ for every $i = 1, \dots, n$.
\end{itemize}

If $X$ is DPT at every $x \in X$, then $X$ is called a {\it DPT space}.
   \index{dpt@@DPT. A metric space $X$ is DPT at $x \in X$}
   \index{dpt@@DPT. A metric space is DPT}

(c) Let $X$ be a metric space and $x \in X$. We shall say that
$X$ has {\em connectivity property 1 at $x$}, (briefly: {\em $X$ is CP1 at $x$}),
if for every $r > 0$ there is $r^* \in (0,r)$
such that for every $x' \in X$ and $r' > 0$:
if $B(x',r') \subs B(x,r^*)$ and $C$ is a connected component
of $B(x,r) \sm B(x',r')$,
then $C \cap (B(x,r) \sm B(x,r^*)) \nnempty$.

If $X$ is CP1 at every $x \in X$, then $X$ is called a {\em CP1 space}.
   \index{CP1 at a point}
   \index{CP1 space}
\end{df}

\begin{tw}[\cite{RY}, Theorem 3.41]\label{t6.2}
Assume that the following facts hold.
\begin{enumerate}
\item[(i)] 
$X$ is a metric space, $G \leq \homeo(X)$,
$\pair{G}{\sigma}$ is a topological group.
The action of $G$ on $X$ is continuous with respect to $\sigma$
and $\pair{G}{\sigma}$ is of the second category.
\item[(ii)] 
$x \in X$ and $\pair{G}{\sigma}$ is affine-like at $x$.
\item[(iii)]
$\Gamma$ is a countably generated modulus of continuity.
\item[(iv)]
$Y$ is a metric space and $\iso{\tau}{X}{Y}$.
\item[(v)]
For every $g \in G$,
$g^{\tau}$ is $\Gamma$-bicontinuous at $\tau(x)$.
\item[(vi)]
$X$ is DPT at $x$ and $Y$ is DPT and CP1 at $\tau(x)$.
\end{enumerate}
{Then} $\tau$ is $\Gamma$-bicontinuous at $x$.
\end{tw}

We need to know that $X$ is DPT, CP1
and that $\lip^{\id{}}(X)$ is affine-like.
The verification of the first two properties is trivial,
and is left to the reader.
We only prove the affine-likeness of $\lip^{\id{}}(X)$.

\begin{lm}\label{l6.3}
Let $X$ be a nonempty open subset of $\bbU$. Then

(a) $X$ is DPT and CP1.

(b) $\lip^{\id{}}(X)$ is affine-like.
\end{lm}

\begin{pf}
(b) Let $u \in X$
and $V \in \nbd^{\lip^{\id{}}(X)}(\id{})$.
We prove that there exists $U \in \nbd^X(u)$
such that for every $x_1,y_1,x_2,y_2 \in U$:
if $d(x_1,y_1) = d(x_2,y_2)$, then there is $g \in V$ such that
$g(x_1) = x_2$ and $g(y_1) = y_2$.
Note that this implies the affine-likeness of
$\lip^{\id{}}(X)$
with $n(x,V,W) = 1$ for every $x,V$ and $W$ (see Definition \ref{d6.1}).

We may assume that $V = B(\id{},s)$.
Choose $K \in (1,e^s)$.
Note that for every $g \in \lip^{\id{}}(X)$,
if $g$ is $K$-bilipschitz,
and $d(x,g(x)) < s$ for every $x \in X$,
then $g \in V$.
Let $r_0 > 0 $ be such that $B^{\bbU}(u,r_0) \subs X$.
Let $a = \min(\frac{r_0}{4},\frac{s}{K + 1})$
and $b = \frac{(K - 1) \cdot a}{(2K - 1) \cdot (K + 1)}$.
Let $u_0$ be such that $d(u_0,u) = a$ and
$B = B^{\bbU}(u_0,2a)$ and $U = B^{\bbU}(u,b)$.

Note that $B \subs B^{\bbU}(u,r_0) \subs X$, because $3a < r_0$.
Note also that $U \subs B$
and that $d(U,\bbU \sm B) = a - b$.
This is so,
since the radius of $B$ is $2a$ and $U$ is a ball whose center
has distance $a$ from the center of $B$ and whose radius is $b$.
Let $f$ be a one-to-one function such that
$\dom(f)$ and $\rng(f)$ are finite subsets of $U$.
We estimate from above
$$
m(f) := \max\left(\bigsetof{\frac{d(x,f(x))}{2a - d(x,u_0)}}
{x \in \dom(f)}\right).
$$
Clearly, $d(x,u_0) \leq d(u_0,u) + d(u,x) < a + b$.
So $2a - d(x,u_0) > a - b$.
Hence 
$m(f) < \frac{2b}{a - b}$.
Define $N := \frac{a - b}{2b}$.
Since $m(f^{-1}) < \frac{2b}{a - b}$, it follows that
$f$ is $N$-bigood.
Suppose that $f$ is an isometry and let
$g = f \cup \sn{\pair{u_0}{u_0}}$.
We show that $g$
is $K$-bilipschitz.
Since the assumptions about $f$ and $f^{-1}$ are the same,
it suffices to check that
$\frac{d(u_0,f(x))}{d(u_0,x)} \leq K$
for every $x \in \dom(f)$.
Note that
$\frac{d(u_0,f(x))}{d(u_0,x)} \leq
\frac{a + b}{a - b}$.
Clearly,
$b = \frac{(K - 1) \cdot a}{(2K - 1) \cdot (K + 1)} \leq
\frac{(K - 1) \cdot a}{K + 1}$.
So
$$
\frac{d(u_0,f(x))}{d(u_0,x)} \leq
\frac{a + b}{a - b} \leq
\frac{a + \frac{(K - 1) \cdot a}{K + 1}}
{a - \frac{(K - 1) \cdot a}{K + 1}} =
\frac{K + 1 + K - 1}{K + 1 - (K - 1)} = K.
$$
We have shown that $g$
is $K$-bilipschitz.
We have also shown that $f$ is $N$-bigood
and so $g$ is $N$-bigood.

We shall now apply the Bilipschitz Extension Theorem to $g$.
For this we still need to show that
$N \geq \frac{K^2}{K - 1}$. Indeed, we have
$$
N = \frac{a - b}{2b} =
\frac{a}{2b} - \frac12 =
\frac{(2K - 1) \cdot (K + 1)}{2(K - 1)} - \frac12 =
\frac{2K^2 + K - 1 - (K - 1)}{2(K - 1)} = \frac{K^2}{K - 1}.
$$
By the Bilipschitz Extension Theorem,
there is $\til h \in \stablep{\homeo(\bbU)}{B}$
such that $g \subs \til h$
and $\til h$ is $K$-bilipschitz.
Hence $h := \til h \rest X \in \lip^{\id{}}(X)$.
We show that for every $z \in \bbU$,
$d(\til h(z),z) < s$.
If $z \not\in B$, then $d(\til h(z),z) = 0$.
Suppose that $z \in B$.
Then
$$d(\til h(z),z) \leq d(z,u_0)  + d(u_0,\til h(z)) <
a + Ka \leq (K + 1) \cdot \frac{s}{K + 1} = s.$$
It follows that $d_S(\til h,\id{}) < s$. Also, $\lipc(\til h) \leq K$, so $d_L(\til h,\id{}) \leq \log(K) < s$. Hence $\hat d(\til h,\id{}) < s$. That is, $h \in B(\id{},s) = V$.
\end{pf}

\begin{wn}\label{c6.4}
(a) Let $X,Y$ be nonempty open subsets of $\bbU$,
$\Gamma$ and $\Delta$ be countably generated \mc-semigroups
and
$\iso{\phi}
{\homeolc{\Gamma}(X)}{\homeolc{\Delta}(Y)}$.
Then there is $\iso{\tau}{X}{Y}$
such that $\tau$ is locally $\Gamma$-bicontinuous,
$\tau$ is locally $\Delta$-bicontinuous
and $\phi(g) = \tau \cmp g \cmp \tau^{-1}$
for every $g \in \homeolc{\Gamma}(X)$.

(b) If $X = Y = \bbU$ in (a), then $\Gamma = \Delta$.
\end{wn}

% CHECK:
\begin{pf}
(a) By Corollary~\ref{c5.10}(a),
$\pair{X}{\homeolc{\Gamma}(X)},
\pair{Y}{\homeolc{\Delta}(Y)} \in K_F$.
So by Theorem~\ref{t5.1}, there is $\iso{\tau}{X}{Y}$
such that $\tau$ induces $\phi$.

Let $G = \lip^{\id{}}(X)$.
As $G$ is a subgroup of $\lip(X)$, it inherits the topology
defined on $\lip(X)$ in Theorem~\ref{t4.2}(c).
Denote this topology on $G$ by $\sigma$.
Since $G \leq \homeolc{\Gamma}(X)$,
it follows that $G^{\tau} \subs \homeolc{\Delta}(Y)$.
We shall show that Theorem~\ref{t6.2} can be applied to
$X,Y,G,\sigma,\tau$ and $\Delta$.

We verify that Clause (i) in Theorem~\ref{t6.2} is fulfilled.
By Theorem~\ref{t4.2}(a) and (b),
$\pair{G}{\sigma}$ is a topological group acting continuously on $X$.
It is easy to see that $G$ is a closed subset of $\lip(X)$.
By Theorem~\ref{t4.2}(c),
$\pair{G}{\sigma}$ is of the second category.

Clause (ii) follows from Lemma~\ref{l6.3}(b),
and Clause~(vi) follows from Lemma~\ref{l6.3}(a).
The remaining requirements of Theorem~\ref{t6.2} hold automatically.
It follows from Theorem~\ref{t6.2} that $\tau$ is locally
$\Delta$-bicontinuous.
Applying the same argument to $\tau^{-1}$ we conclude that
$\tau^{-1}$ is locally $\Gamma$-bicontinuous.
So $\tau$ is locally $\Gamma$-bicontinuous.

(b) It follows from (a) that
$\homeolc{\Gamma}(\bbU) = \homeolc{\Delta}(\bbU)$.
We show that this implies that $\Gamma = \Delta$.
Suppose that $\gamma \in \Gamma \sm \Delta$,
and we shall show that
$\homeolc{\Gamma}(\bbU) \sm \homeolc{\Delta}(\bbU) \neq
\emptyset$.
Let $\setof{\delta_i}{i \in \Nat}$ be a generating set for $\Delta$.
For every $i \in \Nat$ let
$\sett{t_{n,i}}{n \in \Nat} \subs (0,\infty)$
be a sequence converging to $0$ so that 
$\gamma(t_{n,i}) > \delta_i(t_{n,i})$
for every $i$ and~$n$.
It is easy to choose one-to-one sequences
$\sett{x_j}{j \in \Nat}$ and $\sett{y_j}{j \in \Nat}$
so that
\begin{itemize}
\item[(1)] $\lim_{j \to\infty} x_j = x$,
\item[(2)] $d(y_j,x) = \gamma(d(x_j,x))$ for every $j \in \Nat$,
\item[(3)] for every $i \in \Nat$ the set
$\setof{n}{\mbox{there is $j$ such that } d(x_j,x) = t_{n,i}}$
is infinite,
\item[(4)] the function $f$ defined by $x \mapsto x$ and $x_j \mapsto y_j$, $j \in \Nat$,
is $(2 \gamma)$-bicontinuous.
\end{itemize}
Since $\dom(f)$ and $\rng(f)$ are convergent sequences,
they are totally bounded.
Hence by Theorem~\ref{t3.4} there is $g \in \homeo(\bbU)$
such that $g$ is $(2\gamma)$-bicontinuous
and $g \sups f$.
It follows that $g \in \homeolc{\Gamma}(\bbU)$.
However, for every $i \in \Nat$,
$f$ is not $\delta_i$-continuous at $x$.
So $g$ is not $\delta_i$-continuous at $x$.
This means that $g \not\in \homeolc{\Delta}(\bbU)$.
\end{pf}

\section*{Acknowledgements}

The first author would like to thank the Center for Advanced Studies in Mathematics at Ben Gurion University of the Negev, for supporting his visits when this work originated.

%\printindex

\end{document}